\documentclass[12pt, A4paper]{article}
\usepackage{latexsym,amssymb,amsmath}
\usepackage{amssymb,amsthm,upref,amscd}
\usepackage{color}
\usepackage[T1]{fontenc}
\usepackage{indentfirst}
\begin{document}
\baselineskip=14pt
\numberwithin{equation}{section}

\newtheorem{thm}{Theorem}[section]
\newtheorem{lem}[thm]{Lemma}
\newtheorem{cor}[thm]{Corollary}
\newtheorem{Prop}[thm]{Proposition}
\newtheorem{Def}[thm]{Definition}
\newtheorem{Rem}[thm]{Remark}
\newtheorem{Ex}[thm]{Example}

\newcommand{\A}{\mathbb{A}}
\newcommand{\B}{\mathbb{B}}
\newcommand{\C}{\mathbb{C}}
\newcommand{\D}{\mathbb{D}}
\newcommand{\E}{\mathbb{E}}
\newcommand{\F}{\mathbb{F}}
\newcommand{\G}{\mathbb{G}}
\newcommand{\I}{\mathbb{I}}
\newcommand{\J}{\mathbb{J}}
\newcommand{\K}{\mathbb{K}}
\newcommand{\M}{\mathbb{M}}
\newcommand{\N}{\mathbb{N}}
\newcommand{\Q}{\mathbb{Q}}
\newcommand{\R}{\mathbb{R}}
\newcommand{\T}{\mathbb{T}}
\newcommand{\U}{\mathbb{U}}
\newcommand{\V}{\mathbb{V}}
\newcommand{\W}{\mathbb{W}}
\newcommand{\X}{\mathbb{X}}
\newcommand{\Y}{\mathbb{Y}}
\newcommand{\Z}{\mathbb{Z}}
\newcommand\ca{\mathcal{A}}
\newcommand\cb{\mathcal{B}}
\newcommand\cc{\mathcal{C}}
\newcommand\cd{\mathcal{D}}
\newcommand\ce{\mathcal{E}}
\newcommand\cf{\mathcal{F}}
\newcommand\cg{\mathcal{G}}
\newcommand\ch{\mathcal{H}}
\newcommand\ci{\mathcal{I}}
\newcommand\cj{\mathcal{J}}
\newcommand\ck{\mathcal{K}}
\newcommand\cl{\mathcal{L}}
\newcommand\cm{\mathcal{M}}
\newcommand\cn{\mathcal{N}}
\newcommand\co{\mathcal{O}}
\newcommand\cp{\mathcal{P}}
\newcommand\cq{\mathcal{Q}}
\newcommand\rr{\mathcal{R}}
\newcommand\cs{\mathcal{S}}
\newcommand\ct{\mathcal{T}}
\newcommand\cu{\mathcal{U}}
\newcommand\cv{\mathcal{V}}
\newcommand\cw{\mathcal{W}}
\newcommand\cx{\mathcal{X}}
\newcommand\ocd{\overline{\cd}}

\def\c{\centerline}
\def\ov{\overline}
\def\emp {\emptyset}
\def\pa {\partial}
\def\bl{\setminus}
\def\op{\oplus}
\def\sbt{\subset}
\def\un{\underline}
\def\al {\alpha}
\def\bt {\beta}
\def\de {\delta}
\def\Ga {\Gamma}
\def\ga {\gamma}
\def\lm {\lambda}
\def\Lam {\Lambda}
\def\om {\omega}
\def\Om {\Omega}
\def\sa {\sigma}
\def\vr {\varepsilon}
\def\va {\varphi}

\title{\bf \Large Fine bounds for best constants of fractional subcritical Sobolev embeddings and applications to nonlocal PDEs}

\author{Daniele Cassani and Lele Du\\
\small Dip. di Scienza e Alta Tecnologia, Universit\`{a} degli Studi dell'Insubria,\\
\small and \\
\small Riemann International School of Mathematics,\\
\small villa Toeplitz via G.B. Vico 46 - 21100 Varese, Italy.\\
\small \texttt{daniele.cassani@uninsubria.it}, \texttt{dldu@uninsubria.it}
}

\date{\today}

\maketitle

\begin{abstract}
\noindent We establish fine bounds for best constants of the fractional subcritical Sobolev embeddings
\begin{align*}
W_{0}^{s,p}\left(\Omega\right)\hookrightarrow L^{q}\left(\Omega\right),
\end{align*}
where $N\geq1$, $0<s<1$, $p=1,2$, $1\leq q<p_{s}^{\ast}=\frac{Np}{N-sp}$ and $\Omega\subset\mathbb{R}^{N}$ is a bounded smooth domain or the whole space $\mathbb{R}^{N}$. Our results cover the borderline case $p=1$, the Hilbert case $p=2$, $N>2s$ and the so-called Sobolev limiting case $N=1$, $s=\frac{1}{2}$ and $p=2$, where a sharp asymptotic estimate is given by means of a limiting procedure. We apply the obtained results to prove existence and non-existence of solutions for a wide class of nonlocal partial differential equations.
 \vspace{0.3cm}


\vspace{0.3cm}

\end{abstract}

\section{Introduction}
In the study of Partial Differential Equations via a variational approach, the first step is to find a proper function space setting in which the energy functional is well defined and smooth enough to set up equivalence between (weak) solutions to Euler-Lagrange equations and critical points of related functionals: this often yields to consider Sobolev's spaces. Once the underlying functional setting is available, from one side nonlinearities which can be handled, in terms of growth at infinity and near zero, are classified by embedding properties of the function space into other spaces, typically Lebesgue spaces, interpolation spaces between Lebesgue spaces, such as Lorentz spaces, and more general rearrangement invariant spaces \cite{GT}. On the other side, optimal constants involved in integral inequalities responsible of the function space embeddings, turn out to be a kind of ``DNA'' building blocks which can be used to describe qualitative as well as quantitative compactness features. Moreover, borderline cases in Sobolev embeddings have deep connections with geometric measure theory \cite{Talenti} and conformal geometry \cite{Aubin}, and for this reason intensively studied by several authors during the last fifty years.

Nevertheless, there is an aspect which just in recent years has attracted attention, which concerns non-borderline cases where the explicit knowledge of optimal Sobolev's constants is out of reach though still connected to compactness properties of PDEs, as we are going to develop here.

\medskip

\noindent The classical Sobolev constant $\mathcal{S}_{1,p_{1}^{\ast}}$, which was obtained by Aubin \cite{Aubin} and Talenti \cite{Talenti} and explicitly given by
$$
\mathcal{S}_{1,p_{1}^{\ast}}=\pi^{\frac{p}{2}}N\left(\frac{N-p}{p-1}\right)^{p-1}
\left[\frac{\Gamma\left(\frac{N}{p}\right)\Gamma\left(1+N-\frac{N}{p}\right)}{\Gamma\left(\frac{N}{2}+1\right)\Gamma(N)}\right]^{\frac{p}{N}},
$$
appears as the optimal constant of the critical Sobolev embedding
$$
D^{1,p}\left(\mathbb{R}^{N}\right)\hookrightarrow L^{p_{1}^{\ast}}\left(\mathbb{R}^{N}\right),\quad u\in D^{1,p}\left(\mathbb{R}^{N}\right)\backslash\left\{0\right\},
$$
where $N\geq2$, $1<p<N$ and for the critical Sobolev exponent $p_{1}^{\ast}=\frac{Np}{N-p}$, in the sense that $\mathcal{S}_{1,p_{1}^{\ast}}$ is the best possible constant of the inequality
\begin{align}\label{dc1}
C\left\|u\right\|^{p}_{L^{p_{1}^{\ast}}\left(\mathbb{R}^{N}\right)}\leq\left\|\nabla u\right\|^{p}_{L^{p}\left(\mathbb{R}^{N}\right)},\quad u\in D^{1,p}\left(\mathbb{R}^{N}\right)\backslash\left\{0\right\}.
\end{align}
The fact that $\mathcal{S}_{1,p_{1}^{\ast}}$ is explicitly known is due to the invariance property by the group action of dilations and scaling of \eqref{dc1}. As a consequence, the Sobolev constant $\mathcal{S}_{1,p_{1}^{\ast}}$ retains important informations in studying the lack of compactness in nonlinear problems. In particular, the energy levels at which energy functionals fail to satisfy the Palais-Smale condition are quantized in terms of multiples of this constant \cite{Struwe}. Furthermore, this is related to the role of threshold of the Sobolev critical exponent for the existence and nonexistence of solutions to nonlinear PDEs. In fact, the infinitesimal generator of the group invariance yields the so-called Pohozaev-type identities and in turn nonexistence results in fairly smooth domains. This mathematical evidences reflect geometric as well as physical phenomena; see \cite{Berestycki-Lions,Brezis-Coron-Nirenberg,Brezis-Nirenberg,Cerami-Fortunato-Struwe,Lions2,Pohozaev,Pucci-Serrin,Struwe} and references therein.

When $p=1$, it is well-known from geometric measure theory that the Sobolev constant $\mathcal{S}_{1,p_{1}^{\ast}}$ is equal to the isoperimetric constant
$$
\mathcal{S}_{1,1_{1}^{\ast}}=N\omega_{N}^{\frac{1}{N}},
$$
where $\omega_{N}$ is the volume of unit ball, namely
$$
\omega_{N}=\frac{\pi^{\frac{N}{2}}}{\Gamma\left(\frac{N}{2}+1\right)},
$$
see Federer-Fleming \cite{Federer-Fleming}, Fleming-Rishel \cite{Fleming-Rishel} and Maz'ya \cite{Mazya1,Mazya2}.

\medskip

For the subcritical Sobolev embedding
\begin{align*}
W_{0}^{1,p}\left(\Omega\right)\hookrightarrow L^{q}\left(\Omega\right),
\end{align*}
where $N\geq2$, $1\leq p\leq N$, $1\leq q<p_{1}^{\ast}$ and $\Omega\subset\mathbb{R}^{N}$ is a bounded smooth domain or the whole space $\mathbb{R}^{N}$, there still exists optimal constants $S_{1,q}$ in the following inequalities
\begin{align*}
\left\{
\begin{array}{lcl}
\displaystyle C\left\|u\right\|^{p}_{L^{q}\left(\Omega\right)}\leq\left\|\nabla u\right\|^{p}_{L^{p}\left(\Omega\right)},&&u\in W_{0}^{1,p}\left(\Omega\right)\backslash\left\{0\right\}; \\
\\
\displaystyle C\left\|u\right\|^{p}_{L^{q}\left(\mathbb{R}^{N}\right)}\leq\left\|\nabla u\right\|^{p}_{L^{p}\left(\mathbb{R}^{N}\right)}+\left\|u\right\|^{p}_{L^{p}\left(\mathbb{R}^{N}\right)},&&u\in W^{1,p}\left(\mathbb{R}^{N}\right)\backslash\left\{0\right\}.
\end{array}
\right.
\end{align*}
The attainability of the Sobolev constants $S_{1,q}$ is well-known in the literature, whereas there is no hope in general to obtain their explicit value. It is a general fact, the absence of explicit solutions to general nonlinear equations. However, recent applications assume some sharp growth conditions which involve the explicit knowledge of the Sobolev constant $S_{1,q}$, see \cite{Alves-Cassani-Tarsi-Yang,Alves-Soares,Alves-Souto-Montenegro,Cao,Chen-Zou,Ding-Ni,Lions1,O-Souto} and also \cite{Ferone-Murat,Gazzola-Sperone-Weth} for more applications in different contexts. So that it seems to get consolidating a new method which makes a systematic use of growth conditions which involve the best constants $S_{1,q}$. This motivates to searching for fine bounds for $S_{1,q}$ as first established in \cite{Cassani-Tarsi-Zhang} for the Hilbert case $p=2$ and then extended in \cite{Du} up to the general case $1\leq p\leq N$. 

\medskip

Here we are concerned with the fractional Sobolev embeddings
\begin{align*}
W_{0}^{s,p}\left(\Omega\right)\hookrightarrow L^{q}\left(\Omega\right),
\end{align*}
where $N\geq1$, $0<s<1\leq p\leq\frac{N}{s}$ and $q$ satisfies
\begin{align*}
\left\{
\begin{array}{lcl}
\displaystyle 1\leq q\leq p_{s}^{\ast},&&N>sp,~\Omega~\text{is}~\text{bounded}; \\
\\
\displaystyle 1\leq q<+\infty,&&N=sp,~\Omega~\text{is}~\text{bounded}; \\
\\
\displaystyle p\leq q\leq p_{s}^{\ast},&&N>sp,~\Omega=\mathbb{R}^{N}; \\
\\
\displaystyle p\leq q<+\infty,&&N=sp,~\Omega=\mathbb{R}^{N}
\end{array}
\right.
\end{align*}
with the fractional Sobolev critical exponent
$$
p_{s}^{\ast}=\frac{Np}{N-sp}.
$$

For the fractional critical Sobolev embedding
$$
D^{s,p}\left(\Omega\right)\hookrightarrow L^{p_{s}^{\ast}}\left(\Omega\right),
$$
there exists an optimal constant $S_{s,p_{s}^{\ast}}\left(\Omega\right)$ such that
$$
S_{s,p_{s}^{\ast}}\left(\Omega\right)\left\|u\right\|^{p}_{L^{p^{\ast}_{s}}\left(\Omega\right)}
\leq\left[u\right]_{W^{s,p}\left(\mathbb{R}^{N}\right)}^{p},\quad u\in D^{s,p}\left(\Omega\right)\backslash\left\{0\right\},
$$
where $\left[u\right]_{W^{s,p}\left(\mathbb{R}^{N}\right)}$ is the standard Gagliardo semi-norm, namely
$$
S_{s,p_{s}^{\ast}}\left(\Omega\right)=\inf_{u\in D^{s,p}\left(\Omega\right)\backslash\left\{0\right\}}
\frac{\left[u\right]_{W^{s,p}\left(\Omega\right)}^{p}}{\left\|u\right\|^{p}_{L^{p_{s}^{\ast}}\left(\Omega\right)}}.
$$
The invariance by scaling of the quotient $S_{s,p_{s}^{\ast}}\left(\Omega\right)$ implies that $S_{s,p_{s}^{\ast}}\left(\Omega\right)$ is independent of $\Omega$ and thus
$$
S_{s,p_{s}^{\ast}}\left(\Omega\right)=S_{s,p_{s}^{\ast}}\left(\mathbb{R}^{N}\right)=:\mathcal{S}_{s,p_{s}^{\ast}}.
$$
In the borderline case $p=1$, the fractional isoperimetric constant $\mathcal{S}_{s,1_{s}^{\ast}}$ was given by Brasco-Lindgren-Parini \cite{Brasco-Lindgren-Parini}, namely
$$
\mathcal{S}_{s,1_{s}^{\ast}}=\omega_{N}^{\frac{s-N}{N}}\left[\chi_{B_{1}}\right]_{W^{s,1}\left(\mathbb{R}^{N}\right)},
$$
where $\left[\chi_{B_{1}}\right]_{W^{s,1}\left(\mathbb{R}^{N}\right)}$ is the nonlocal s-perimeter of the unit ball $B_{1}$. More precisely, the explict value of $\mathcal{S}_{s,1_{s}^{\ast}}$ can be computed by the results of Frank-Seiringer \cite{Frank-Seiringer}, namely
$$
\mathcal{S}_{s,1_{s}^{\ast}}=\frac{\omega_{N}^{\frac{s}{N}}N}{N-s}A\left(N,s\right),
$$
where $A\left(N,s\right)$ is the sharp constant of the fractional Hardy-Sobolev inequality
$$
A\left(N,s\right)=2\int^{1}_{0}r^{s-1}\left(1-r^{N-s}\right)A_{N,s}\left(r\right)dr
$$
and where
\begin{align*}
A_{N,s}\left(r\right)=\left\{
\begin{array}{lcl}
\displaystyle\left(N-1\right)\omega_{N-1}\int^{1}_{-1}
\frac{\left(1-t^{2}\right)^{\frac{N-3}{2}}}{\left(1-2rt+r^{2}\right)^{\frac{N+s}{2}}}dt,&&N\geq2, \\
\\
\displaystyle\frac{1}{\left(1-r\right)^{1+s}}+\frac{1}{\left(1+r\right)^{1+s}},  &&N=1.
\end{array}
\right.
\end{align*}
The fractional isoperimetric constant $\mathcal{S}_{s,1_{s}^{\ast}}$ is achieved by a scalar multiple of the characteristic function of a ball in $\mathbb{R}^{N}$. In the Hilbert case $p=2$, Lieb \cite{Lieb} computed the Sobolev constant
$$
\mathcal{S}_{s,2_{s}^{\ast}}=\frac{2\pi^{\frac{N}{2}+s}}{s(1-s)}
\left[\frac{\Gamma\left(2-s\right)}{\Gamma\left(\frac{N}{2}-s\right)}\right]
\left[\frac{\Gamma\left(\frac{N}{2}\right)}{\Gamma\left(N\right)}\right]^{\frac{2s}{N}}.
$$
In analogy to $\mathcal{S}_{1,2_{1}^{\ast}}$,
the Sobolev constant $\mathcal{S}_{s,2_{s}^{\ast}}$ appears as a key ingredient in studying the lack of compactness in fractional problems as developed by Servadei-Valdinoci \cite{Servadei-Valdinoci}. The extremal functions of $\mathcal{S}_{s,2_{s}^{\ast}}$ in $D^{s,2}\left(\mathbb{R}^{N}\right)$ were obtained by Lieb \cite{Lieb} and up to translation and dilation, given by
$$
U_{s,2_{s}^{\ast}}\left(x\right)=\left(\frac{1}{1+\left|x\right|^{2}}\right)^{\frac{N-2s}{2}},\quad x\in \mathbb{R}^{N},
$$
whereas $\mathcal{S}_{s,2_{s}^{\ast}}$ has no positive minimizer on any star-shaped domain $\Omega\neq\mathbb{R}^{N}$ by the validity of a fractional Pohozaev-type identity obtained by Ros Oton-Serra \cite{Ros-Oton-Serra}.

However, nothing is known for $\mathcal{S}_{s,p_{s}^{\ast}}$ when $p\in\left(1,2\right)\cup\left(2,+\infty\right)$. Indeed, when $p\neq 2$, the Sobolev space $W^{s,p}\left(\mathbb{R}^{N}\right)$ and the Bessel potential spaces $H^{s,p}\left(\mathbb{R}^{N}\right)$ are no longer equivalent,
and this is a major difficulty to compute $\mathcal{S}_{s,p^{\ast}_{s}}$ by exploiting the sharp Hardy-Littlewood-Sobolev inequality.
A lower bound for $\mathcal{S}_{s,p_{s}^{\ast}}$ was given by Maz'ya-Shaposhnikova \cite{Mazya-Shaposhnikova}
$$
\mathcal{S}_{s,p_{s}^{\ast}}\geq\frac{\omega_{N}N\left(N-sp\right)^{p-1}}{2^{\left(N+1\right)\left(N+2\right)}s
\left(1-s\right)p^{p+2}\left(N+2p\right)^{3p}}.
$$

For the fractional subcritical Sobolev embedding
$$
W_{0}^{s,p}\left(\Omega\right)\hookrightarrow L^{q}\left(\Omega\right),
$$
there exists optimal constants $S_{s,q}$ in the following inequalities
\begin{align*}
\left\{
\begin{array}{lcl}
\displaystyle C\left\|u\right\|^{p}_{L^{q}\left(\Omega\right)}
\leq\left[u\right]_{W^{s,p}\left(\mathbb{R}^{N}\right)}^{p},&&u\in W_{0}^{s,p}\left(\Omega\right)\backslash\left\{0\right\}; \\
\\
\displaystyle C\left\|u\right\|^{p}_{L^{q}\left(\mathbb{R}^{N}\right)}
\leq\left[u\right]_{W^{s,p}\left(\mathbb{R}^{N}\right)}^{p}
+\left\|u\right\|^{p}_{L^{p}\left(\mathbb{R}^{N}\right)},&&u\in W^{s,p}\left(\mathbb{R}^{N}\right)\backslash\left\{0\right\},
\end{array}
\right.
\end{align*}
namely
\begin{align*}
\left\{
\begin{array}{lcl}
\displaystyle S_{s,q}\left(\Omega\right)=\inf_{u\in W^{s,p}_{0}\left(\Omega\right)\backslash\left\{0\right\}}
\frac{\left[u\right]_{W^{s,p}\left(\mathbb{R}^{N}\right)}^{p}}
{\left\|u\right\|^{p}_{L^{q}\left(\Omega\right)}}; \\
\\
\displaystyle S_{s,q}\left(\mathbb{R}^{N}\right)=\inf_{u\in W^{s,p}\left(\mathbb{R}^{N}\right)\backslash\left\{0\right\}}
\frac{\left[u\right]_{W^{s,p}\left(\mathbb{R}^{N}\right)}^{p}
+\left\|u\right\|^{p}_{L^{p}\left(\mathbb{R}^{N}\right)}}{\left\|u\right\|^{p}_{L^{q}\left(\mathbb{R}^{N}\right)}}.
\end{array}
\right.
\end{align*}
The action of the dilation group $u=u\left(\lambda x\right)$ for the quotient $S_{s,q}\left(\Omega\right)$ yields
\begin{align}\label{A}
S_{s,q}\left(\Omega\right)=
\left\{
\begin{array}{lcl}
\displaystyle \lambda^{Np\left(\frac{1}{p^{\ast}_{s}} -\frac{1}{q}\right)}S_{s,q}\left(\Omega_{\lambda}\right),&&N>ps; \\
\\
\displaystyle \lambda^{-\frac{Np}{q}}S_{s,q}\left(\Omega_{\lambda}\right),&&N=ps,
\end{array}
\right.
\end{align}
which means that $S_{1,q}\left(\Omega\right)$ strictly depends on the domain $\Omega$ when $1\leq q<p_{s}^{\ast}$.

In particular, when $p=2$, in order to be consistent with the definition of the Sobolev constant $\mathcal{S}_{1,2_{1}^{\ast}}$, one replaces the Gagliardo semi-norm $\left[u\right]_{W^{s,2}\left(\mathbb{R}^{N}\right)}$ by an equivalent $L^{2}$ norm of the fractional Laplace operator, hence one can define the optimal constant of the following fractional Sobolev inequality:
$$
S_{s}\left(\Omega\right)\left\|u\right\|^{2}_{L^{2^{\ast}_{s}}\left(\Omega\right)}
\leq\left\|\left(-\Delta\right)^{\frac{s}{2}}u\right\|^{2}_{L^{2}\left(\mathbb{R}^{N}\right)},\quad u\in D^{s,2}\left(\Omega\right)\backslash\left\{0\right\},
$$
namely
$$
S_{s}\left(\Omega\right)=\inf_{u\in D^{s,2}\left(\Omega\right)\backslash\left\{0\right\}}
\frac{\left\|\left(-\Delta\right)^{\frac{s}{2}}u\right\|^{2}_{L^{2}\left(\mathbb{R}^{N}\right)}}
{\left\|u\right\|^{2}_{L^{2^{\ast}_{s}}\left(\Omega\right)}}.
$$
After applying the identity
\begin{align}\label{B}
\left[u\right]_{W^{s,2}\left(\mathbb{R}^{N}\right)}^{2}=\frac{2}{B\left(N,s\right)}
\left\|\left(-\Delta\right)^{\frac{s}{2}}u\right\|^{2}_{L^{2}\left(\mathbb{R}^{N}\right)},
\end{align}
where
$$
B\left(N,s\right)=\frac{2^{2s}s}{\pi^{\frac{N}{2}}}\left[\frac{\Gamma\left(\frac{N}{2}+s\right)}{\Gamma\left(1-s\right)}\right],
$$
the Sobolev constant $\mathcal{S}_{s}:=S_{s}\left(\Omega\right)=S_{s}\left(\mathbb{R}^{N}\right)$ is given by
$$
\mathcal{S}_{s}=2^{2s}\pi^{s}
\left[\frac{\Gamma\left(\frac{N}{2}+s\right)}{\Gamma\left(\frac{N}{2}-s\right)}\right]
\left[\frac{\Gamma\left(\frac{N}{2}\right)}{\Gamma(N)}\right]^{\frac{2s}{N}}\ .
$$
Let us mention that one can also apply the dual property of Hardy-Littlewood-Sobolev inequality like Cotsiolis-Tavoularis \cite{Cotsiolis-Tavoularis} to get the same value. Notice that $\mathcal{S}_{s}\rightarrow\mathcal{S}_{1,2_{1}^{\ast}}$ as $s\rightarrow1^{-}$, hence $\mathcal{S}_{s}$ can be regarded as a generalization of the Sobolev constant $\mathcal{S}_{1,2_{1}^{\ast}}$.

Moreover, for the fractional subcritical Sobolev embedding
$$
H_{0}^{s}\left(\Omega\right)\hookrightarrow L^{q}\left(\Omega\right),
$$
we also replace the definition of the optimal constant $S_{s,q}\left(\Omega\right)$ in the following inequality
\begin{align*}
\left\{
\begin{array}{lcl}
\displaystyle C\left\|u\right\|^{2}_{L^{q}\left(\Omega\right)}
\leq\left\|\left(-\Delta\right)^{\frac{s}{2}}u\right\|^{2}_{L^{2}\left(\mathbb{R}^{N}\right)},&&u\in H_{0}^{s}\left(\Omega\right)\backslash\left\{0\right\}; \\
\\
\displaystyle C\left\|u\right\|^{2}_{L^{q}\left(\mathbb{R}^{N}\right)}
\leq\left\|\left(-\Delta\right)^{\frac{s}{2}}u\right\|^{2}_{L^{2}\left(\mathbb{R}^{N}\right)}
+\left\|u\right\|^{2}_{L^{2}\left(\mathbb{R}^{N}\right)},&&u\in H^{s}\left(\mathbb{R}^{N}\right)\backslash\left\{0\right\},
\end{array}
\right.
\end{align*}
namely
\begin{align*}
\left\{
\begin{array}{lcl}
\displaystyle S_{s,q}\left(\Omega\right)=\inf_{u\in H^{s}_{0}\left(\Omega\right)\backslash\left\{0\right\}}
\frac{\left\|\left(-\Delta\right)^{\frac{s}{2}}u\right\|^{2}_{L^{2}\left(\mathbb{R}^{N}\right)}}
{\left\|u\right\|^{2}_{L^{q}\left(\Omega\right)}}; \\
\\
\displaystyle S_{s,q}\left(\mathbb{R}^{N}\right)=\inf_{u\in H^{s}\left(\mathbb{R}^{N}\right)\backslash\left\{0\right\}}
\frac{\left\|\left(-\Delta\right)^{\frac{s}{2}}u\right\|^{2}_{L^{2}\left(\mathbb{R}^{N}\right)}
+\left\|u\right\|^{2}_{L^{2}\left(\mathbb{R}^{N}\right)}}{\left\|u\right\|^{2}_{L^{q}\left(\mathbb{R}^{N}\right)}}.
\end{array}
\right.
\end{align*}
The Sobolev constant $S_{s,q}\left(\Omega\right)$ is always achieved by means of the compact embedding
$$
H_{0}^{s}\left(\Omega\right)\hookrightarrow L^{q}\left(\Omega\right),\quad 1\leq q<2_{s}^{\ast}
$$
and $S_{s,q}\left(\mathbb{R}^{N}\right)$ is achieved when $2<q<2_{s}^{\ast}$ by the existence results of Frank-Lenzmann \cite{Frank-Lenzmann} for $N=1$ and Dipierro-Palatucci-Valdinoci \cite{Dipierro-Palatucci-Valdinoci} for $N\geq2$, whereas $S_{s,2}\left(\mathbb{R}^{N}\right)$ and $S_{s,2_{s}^{\ast}}\left(\mathbb{R}^{N}\right)$ are never achieved thanks to the fractional Pohozaev-type identity established by Chang-Wang \cite{Chang-Wang} in the whole $\mathbb{R}^{N}$.

The study of quantitative aspects of fractional Sobolev constants is not only interesting from the theoretical point of view. In fact, the Sobolev constant $\mathcal{S}_{s,p_{s}^{\ast}}$, as in the integer case $s=1$, plays an important role in compactness issues and the fractional critical Sobolev exponent $p_{s}^{\ast}$ yields the sharp threshold for the existence and nonexistence of solutions to nonlocal PDEs.

Likewise classical problems \cite{Alves-Cassani-Tarsi-Yang,Alves-Soares,Alves-Souto-Montenegro,Cao,Chen-Zou,Ding-Ni,Lions1,O-Souto}, so far there are plenty of applications \cite{Alves-Figueiredo-Siciliano,Alves-Miyagaki,Alves-O-Miyagaki,Ambrosio,Zhen-He-Xu-Yang} which assume sharp growth conditions involving the explicit knowledge of the Sobolev constants $S_{s,q}$, which turn out to be crucial to determine the existence and nonexistence of solutions to partial differential equations. Those approaches essentially extend the perturbation technique of Br\'{e}zis-Nirenberg \cite{Brezis-Nirenberg} in which a prescribed asymptotic behavior near zero is assumed. Hence, looking for possibly sharp bounds of Sobolev's constants makes such sufficient conditions effective, both from the theoretical point of view and that of applications.

\subsection*{Main results}
\subsubsection*{Bounds for $S_{s,q}\left(\Omega\right)$ and $S_{s,q}\left(\mathbb{R}^{N}\right)$}

Let $N\geq1$, $0<s<1\leq p\leq\frac{N}{s}$ and $\Omega$ be a bounded smooth domain $\Omega\subset\mathbb{R}^{N}$. We denote the largest radius of $\Omega$ by
$$
R_{\Omega}=\sup\left\{R:B_{R}\left(x\right)\subset\Omega,x\in\Omega\right\}
$$
and $\mathrm{B}\left(x,y\right)$ is the Beta function. Let us begin with the case $p=1$.

\begin{thm}\label{thm1}
Let $p=1$ and $1\leq q<1_{s}^{\ast}$. The following hold:
\begin{itemize}
\item [$\left(1\right)$] if $1\leq q<1_{s}^{\ast}$, then
\begin{align*}
\mathcal{S}_{s,1_{s}^{\ast}}\left|\Omega\right|^{\frac{1}{1^{\ast}_{s}}-\frac{1}{q}}
\leq S_{s,q}\left(\Omega\right)\leq\mathcal{S}_{s,1_{s}^{\ast}}\left|B_{R_{\Omega}}\right|^{\frac{1}{1^{\ast}_{s}}-\frac{1}{q}}.
\end{align*}
\item [$\left(2\right)$] $S_{s,1}\left(\mathbb{R}^{N}\right)=1$;
\item [$\left(3\right)$] if $1<q<1_{s}^{\ast}$, then
\begin{align*}
&\left[\frac{N}{s}\left(\frac{1}{q}-\frac{1}{1_{s}^{\ast}}\right)\right]^{\frac{N}{s}\left(\frac{1}{1_{s}^{\ast}}-\frac{1}{q}\right)}
\left[\frac{N}{s}\left(1-\frac{1}{q}\right)\right]
^{\frac{N}{s}\left(\frac{1}{q}-1\right)}
\mathcal{S}_{s,1_{s}^{\ast}}^{\frac{N}{s}\left(1-\frac{1}{q}\right)}\\
&\quad\quad\quad\leq S_{s,q}\left(\mathbb{R}^{N}\right)\leq \frac{s}{N}\left(\frac{1}{q}-\frac{1}{1_{s}^{\ast}}\right)^{\frac{N}{s}\left(\frac{1}{1_{s}^{\ast}}-\frac{1}{q}\right)}
\left(1-\frac{1}{q}\right)^{\frac{N}{s}\left(\frac{1}{q}-1\right)}\mathcal{S}_{s,1_{s}^{\ast}}^{\frac{N}{s}\left(1-\frac{1}{q}\right)}.
\end{align*}
\end{itemize}
\end{thm}

Next we consider the Hilbert case $2=p<\frac{N}{s}$.

\begin{thm}\label{thm2}
Let $2=p<\frac{N}{s}$. The following hold:
\begin{itemize}
\item [$\left(1\right)$] if $1\leq q<2_{s}^{\ast}$, then
\begin{align*}
 \mathcal{S}_{s}\left|\Omega\right|^{2\left(\frac{1}{2^{\ast}_{s}}-\frac{1}{q}\right)}
 \leq S_{s,q}\left(\Omega\right)\leq &\frac{2^{2s+\frac{2}{q}}\left(\omega_{N}N\right)^{1-\frac{2}{q}}}{N+2s}\Gamma^{2}\left(s+1\right)\\
&\cdot\left[\mathrm{B}\left(\frac{N}{2},qs+1\right)\right]^{-\frac{2}{q}}
R_{\Omega}^{2N\left(\frac{1}{2^{\ast}_{s}} -\frac{1}{q}\right)};
\end{align*}
\item [$\left(2\right)$] $S_{s,2}\left(\mathbb{R}^{N}\right)=1$;
\item [$\left(3\right)$] if $2<q<2_{s}^{\ast}$, then
\begin{align*}
&\left[\frac{N}{s}\left(\frac{1}{q}-\frac{1}{2_{s}^{\ast}}\right)\right]^{\frac{N}{s}\left(\frac{1}{2_{s}^{\ast}}-\frac{1}{q}\right)}
\left[\frac{N}{s}\left(\frac{1}{2}-\frac{1}{q}\right)\right]
^{\frac{N}{s}\left(\frac{1}{q}-\frac{1}{2}\right)}
\mathcal{S}_{s}^{\frac{N}{s}\left(\frac{1}{2}-\frac{1}{q}\right)}\\
&\quad\quad\quad\leq S_{s,q}\left(\mathbb{R}^{N}\right)\leq
\omega^{1-\frac{2}{q}}_{N}s
\left[\frac{2^{2s+1-\frac{2s}{N}}\Gamma^{2}\left(s+1\right)}{\left(N+2s\right)
\left(\frac{1}{2}-\frac{1}{q}\right)}\right]^{\frac{N}{s}\left(\frac{1}{2}-\frac{1}{q}\right)}\\
&\quad\quad\quad\quad\quad\quad\quad\quad\quad\quad\cdot\left[N\mathrm{B}\left(\frac{N}{2},qs+1\right)\right]^{-\frac{2}{q}}
\left[\frac{\mathrm{B}\left(\frac{N}{2},2s+1\right)}{\frac{1}{q}-\frac{1}{2^{\ast}_{s}}}\right]
^{\frac{N}{s}\left(\frac{1}{q}-\frac{1}{2^{\ast}_{s}}\right)}.
\end{align*}
\end{itemize}
\end{thm}

In the limiting case $N=2s=1$, a lower bound for $S_{\frac{1}{2},q}\left(\mathbb{R}\right)$ was given by Lieb-Loss \cite{Lieb-Loss}, namely
$$
S_{\frac{1}{2},q}\left(\mathbb{R}\right)\geq\left(q-1\right)^{1-\frac{1}{q}}\left[\frac{q\left(q-2\right)}{2\pi}\right]^{\frac{2}{q}-1}.
$$
When $N=2s=2$, the asymptotic behavior of $S_{1,q}\left(\Omega\right)$ and $S_{1,q}\left(\mathbb{R}^{2}\right)$ were obtained in \cite{Cassani-Tarsi-Zhang,Ren-Wei}, namely
\begin{align}\label{C}
\lim_{q\rightarrow +\infty} qS_{1,q}\left(\Omega\right)=\lim_{q\rightarrow +\infty} qS_{1,q}\left(\mathbb{R}^{2}\right)=8\pi e\ .
\end{align}
Finally, we establish bounds for $S_{\frac{1}{2},q}\left(\Omega\right)$ and $S_{\frac{1}{2},q}\left(\mathbb{R}^{N}\right)$.
\begin{thm}\label{thm3}
Let $s=\frac{1}{2}$, $p=2$ and $N=1$. The following hold:
\begin{itemize}
\item [$\left(1\right)$] if $q\geq 1$, then
\begin{align*}
S_{\frac{1}{2},q}\left(\Omega\right)\leq \frac{2^{1-\frac{2}{q}}\pi e}{q}
R_{\Omega}^{-\frac{2}{q}};
\end{align*}
\item [$\left(2\right)$] $S_{\frac{1}{2},2}\left(\mathbb{R}\right)=1$;
\item [$\left(3\right)$] if $q>2$, then
\begin{align*}
S_{\frac{1}{2},q}\left(\mathbb{R}\right)\leq 2^{1-\frac{4}{q}}\pi^{1-\frac{2}{q}}q\left(q-2\right)^{\frac{4}{q}-2}e^{\frac{q-2}{q}};
\end{align*}
\item [$\left(4\right)$] The asymptotic behavior of $S_{\frac{1}{2},q}\left(\Omega\right)$ and $S_{\frac{1}{2},q}\left(\mathbb{R}^{N}\right)$ is given by
\begin{align*}
\lim_{q\rightarrow +\infty } qS_{\frac{1}{2},q}\left(\Omega\right)=\lim_{q\rightarrow +\infty } qS_{\frac{1}{2},q}\left(\mathbb{R}^{N}\right)=2\pi e.
\end{align*}
\end{itemize}
\end{thm}

\begin{Rem}
\begin{itemize}
\item [1.] In Theorem \ref{thm1}, if $\Omega$ is a ball, then $S_{s,q}\left(\Omega\right)$ is achieved by a scalar multiple of the characteristic function of a ball in $\Omega$ and
\begin{align*}
S_{s,q}\left(\Omega\right)=\mathcal{S}_{s,1_{s}^{\ast}}\left|\Omega\right|^{\frac{1}{1_{s}^{\ast}}-\frac{1}{q}}.
\end{align*}
\item [2.] In Theorem \ref{thm1}, when $q\rightarrow 1_{s}^{\ast}$, we obtain $S_{s,q}\left(\Omega\right)\rightarrow\mathcal{S}_{1,1_{s}^{\ast}}$ and $S_{s,q}\left(\mathbb{R}^{N}\right)\rightarrow\mathcal{S}_{1,1_{s}^{\ast}}$.
\item [3.] When $q\rightarrow 2_{s}^{\ast}$, the lower bound for $S_{s,q}\left(\Omega\right)$ and $S_{s,q}\left(\mathbb{R}^{N}\right)$ in Theorem \ref{thm2} goes to $\mathcal{S}_{s}$.
\end{itemize}
\end{Rem}

\subsubsection*{Applications to nonlocal PDEs}
Let us look for the standing waves solutions $v\left(t,x\right)=e^{i\omega t}u\left(x\right)$ of the following nonlocal nonlinear Schr\"{o}dinger equation
\begin{align*}
iv_{t}=\left(-\Delta\right)^{s} v+\left(V+\omega\right)v-Q\left|v\right|^{q-2}v,
\end{align*}
where $0<s\leq1$, $N\geq2s$, $\omega\in\mathbb{R}$, $Q,V\in \mathcal{C}\left(\mathbb{R}^{N}\right)$ and $q$ satisfies
\begin{align*}
\left\{
\begin{array}{lcl}
\displaystyle 2<q\leq 2_{s}^{\ast},&&N>2s; \\
\\
\displaystyle 2<q<+\infty,&&N=2s,
\end{array}
\right.
\end{align*}
which yields to the following equation
\begin{align}\label{A1}
\left\{
\begin{array}{lcl}
\displaystyle \left(-\Delta\right)^{s} u+Vu=Q\left|u\right|^{q-2}u,\\
\\
\displaystyle u\in H^{s}\left(\mathbb{R}^{N}\right).
\end{array}
\right.
\end{align}
We refer to Laskin \cite{Laskin} for the physical background. When $s=1$, if $V=1$ and $Q$ satisfies
\begin{itemize}
\item [$\left(Q_{1}\right)$] $\displaystyle Q\not\equiv1$;
\item [$\left(Q_{2}\right)$] $\displaystyle Q\geq 1$;
\item [$\left(Q_{3}\right)$] $\displaystyle \lim_{\left|x\right|\rightarrow+\infty} Q\left(x\right)=1$,
\end{itemize}
then Ding-Ni \cite{Ding-Ni} proved \eqref{A1} has a positive solution.
If $Q=1$ and $V$ satisfies
\begin{itemize}
\item [$\left(V_{1}\right)$] $\displaystyle V\not\equiv1$;
\item [$\left(V_{2}\right)$] $\displaystyle 0<V\leq 1$;
\item [$\left(V_{3}\right)$] $\displaystyle \lim_{\left|x\right|\rightarrow+\infty} V\left(x\right)=1$,
\end{itemize}
then Lions \cite{Lions2} proved \eqref{A1} has a positive solution, which is actually a ground state solution. In the results of Ding-Ni \cite{Ding-Ni} and Lions \cite{Lions2}, the Sobolev constant $S_{1,q}\left(\mathbb{R}^{N}\right)$ plays a key role in proving compactness by establishing the existence of a nontrivial limit of a PS sequence which is a solution to the equation. More precisely, let us consider the energy functional
$$
Q_{s,q}\left(u\right)=\frac{1}{2}\int_{\Omega}\left|\left(-\Delta\right)^{\frac{s}{2}} u\right|^{2}dx-\frac{1}{q}\int_{\Omega}\left|u\right|^{q}dx
$$
and the energy level
$$
\beta_{s,q}=\left(\frac{1}{2}-\frac{1}{q}\right)S_{s,q}^{\frac{q}{q-2}}\left(\mathbb{R}^{N}\right).
$$
As $\beta_{1,q}$ yields the first non-compactness level of the energy functional $Q_{1,q}$, the value $\beta_{s,q}$ is the first level of $Q_{s,q}$, where the lack of compactness occurs. We have the following
\begin{thm}\label{thm4}
Let $0<s<1$, $N\geq2s$ and
\begin{align*}
\left\{
\begin{array}{lcl}
\displaystyle 2<q\leq 2_{s}^{\ast},&&N>2s; \\
\\
\displaystyle 2<q<+\infty,&&N=2s=1.
\end{array}
\right.
\end{align*}
\begin{itemize}
\item [$\left(1\right)$] If $V=1$ and $Q$ satisfies $\left(Q_{1}\right)$, $\left(Q_{2}\right)$ and $\left(Q_{3}\right)$, then \eqref{A1} has a positive solution with
$$
\left\|u\right\|^{2}_{H^{s}\left(\mathbb{R}^{N}\right)}<S_{s,q}^{\frac{q}{q-2}}\left(\mathbb{R}^{N}\right);
$$
\item [$\left(2\right)$] If $Q=1$ and $V$ satisfies $\left(V_{1}\right)$, $\left(V_{2}\right)$ and $\left(V_{3}\right)$, then \eqref{A1} has a positive solution with
$$
\left\|u\right\|_{L^{q}\left(\mathbb{R}^{N}\right)}<S_{s,q}^{\frac{1}{q-2}}\left(\mathbb{R}^{N}\right).
$$
\end{itemize}
\end{thm}
Next we consider the following nonlinear and nonlocal scalar field equation:
\begin{align}\label{A2}
\left\{
\begin{array}{lcl}
\displaystyle \left(-\Delta\right)^{s} u+u=f\left(u\right),\\
\\
\displaystyle u\in H^{s}\left(\mathbb{R}^{N}\right),
\end{array}
\right.
\end{align}
where $0<s\leq1$, $N=2s$ and $f\in \mathcal{C}\left(\mathbb{R}\right)$. The equation \eqref{A2} is a special case of the equation
\begin{align}\label{A3}
\left\{
\begin{array}{lcl}
\displaystyle \left(-\Delta\right)^{s} u=g\left(u\right),\\
\\
\displaystyle u\in H^{s}(\mathbb{R}^{N}),
\end{array}
\right.
\end{align}
where $g\left(t\right)=f\left(t\right)-t$. When $s=1$, Berestycki-Lions \cite{Berestycki-Lions} in $N\geq3$ and Berestycki-Gallou\"{e}t-Kavian \cite{Berestycki-Gallouet-Kavian} in $N=2$ studied the constraint minimization problem related to \eqref{A3} with $g$ satisfying subcritical growth conditions and proved the existence of a minimizer which turns out to be a positive ground state solution. After that, Alves-Souto-Montenegro \cite{Alves-Souto-Montenegro} established the existence of a positive ground state solution of \eqref{A2} in $N\geq2$ under the assumptions that $f$ has critical growth conditions and
$$
f\left(t\right)\geq \lambda t^{q-1},\quad t\geq0,
$$
where
\begin{align*}
\lambda\geq\left\{
\begin{array}{lcl}
\displaystyle\left[\frac{N^{\frac{N}{2}}
\left(q-2\right)}{2q\mathcal{S}_{1,2_{1}^{\ast}}^{\frac{N}{2}}\left(N-2\right)^{\frac{N-2}{2}}}\right]
^{\frac{q-2}{2}}S^{\frac{q}{2}}_{1,q}\left(\mathbb{R}^{N}\right),&&N\geq3, \\
\\
\displaystyle\left(\frac{q-2}{q}\right)
^{\frac{q-2}{2}}S^{\frac{q}{2}}_{1,q}\left(\mathbb{R}^{2}\right),&&N=2.
\end{array}
\right.
\end{align*}
This kind of growth conditions closely rely on the lower bound for $S_{1,q}\left(\mathbb{R}^{N}\right)$. When $0<s<1$, the existence of positive ground state solutions of \eqref{A3} was established by Chang-Wang \cite{Chang-Wang} in $N\geq2$ and Alves, Figueiredo and Siciliano \cite{Alves-Figueiredo-Siciliano} extends the results by Alves-Souto-Montenegro \cite{Alves-Souto-Montenegro} to the fractional problem \eqref{A2}, where
\begin{align*}
\lambda\geq\left\{
\begin{array}{lcl}
\displaystyle\left[\frac{N^{\frac{N}{2s}}\left(q-2\right)}{2sq\mathcal{S}_{s}^{\frac{N}{2s}}\left(N-2s\right)^{\frac{N-2s}{2s}}}\right]
^{\frac{q-2}{2}}S^{\frac{q}{2}}_{s,q}(\mathbb{R}^{N}),&&0<s<1,N\geq2, \\
\\
\displaystyle\left(\frac{q-2}{q}\right)
^{\frac{q-2}{2}}S^{\frac{q}{2}}_{s,q}(\mathbb{R}),&&s=\frac{1}{2},N=1.
\end{array}
\right.
\end{align*}
We have the following existence theorem for \eqref{A2}.
 \begin{thm}\label{thm5}
Let $0<s<1$, $N=1$ and $f$ satisfies that for all $t\in\mathbb{R}$, there hold
\begin{itemize}
\item [$\left(f_{1}\right)$] $\displaystyle \lim_{t\rightarrow0^{+}}\frac{f\left(t\right)}{t}=0$;
\item [$\left(f_{2}\right)$] $\displaystyle \left|f(t)\right|\leq C e^{\pi t^{2}}$;
\item [$\left(f_{3}\right)$] $\displaystyle f\left(t\right)\geq \frac{q}{2}S^{\frac{q}{2}}_{\frac{1}{2},q}\left(\mathbb{R}\right)\left|t\right|^{q-2}t$.
\end{itemize}
Then \eqref{A2} has a positive ground state solution.
\end{thm}
Let us stress the fact that growth condition on $f$ as in Theorem \ref{thm5}-$\left(f_{3}\right)$ was just theoretical up to the bound provided in Theorem \ref{thm3} and for instance could not be implemented in numerical applications.

\medskip

Finally, we study the following class of system of strongly coupled nonlocal and nonlinear Schr\"odinger equations
\begin{align}\label{A4}
\left\{
\begin{array}{lcl}
\displaystyle \left(-\Delta\right)^{s} u+u=\left|u\right|^{p-2}u+\lambda v,\\
\\
\displaystyle \left(-\Delta\right)^{s} v+v=\left|v\right|^{q-2}v+\lambda u,\\
\\
\displaystyle u,v\in H^{s}\left(\mathbb{R}^{N}\right),
\end{array}
\right.
\end{align}
where $0<s\leq1$, $N>2s$, $2<p,q\leq2_{s}^*$ and $0<\lambda<1$. When $s=1$, if $N=3$ and $\lambda$ is small, Ambrosetti-Colorado-Ruiz \cite{Ambrosetti-Colorado-Ruiz} proved that there exists multi-bump solitons in \eqref{A4}, provided $2<p=q<2_{1}^*$ and $0<\lambda<1$ and then in Ambrosetti-Cerami-Ruiz \cite{Ambrosetti-Cerami-Ruiz} the authors prove that \eqref{A4} has a positive ground state solution. If $N\geq3$, $2<p,q<2_{1}^*$ and $0<\lambda<1$, Br\'{e}zis-Lieb \cite{Brezis-Lieb} proved \eqref{A4} has a positive ground state solution.

Chen-Zou in \cite{Chen-Zou} established the following
\begin{thm}\label{thm6}
Let $s=1$, $N\geq3$, and $0<\lambda<1$.
\begin{itemize}
\item [$\left(1\right)$] if $2<p<2_{1}^*$ and $q=2_{1}^*$, let
$$
\alpha_{1}=\left[\frac{\mathcal{S}_{1,2_{1}^{\ast}}^{\frac{N}{2}}}{N\left(\frac{1}{2}-\frac{1}{q}\right)
S_{1,q}^{\frac{q}{q-2}}\left(\mathbb{R}^{N}\right)}\right]^{\frac{1}{\frac{q}{q-2}-\frac{N}{2}}},
$$
then there exists $\lambda_{1}\in\left[\sqrt{1-\alpha_{1}},1\right)$ such that
\begin{itemize}
\item [$\left(i\right)$] if $\lambda<\lambda_{1}$, then \eqref{A4} has no ground state solution;
\item [$\left(ii\right)$] if $\lambda>\lambda_{1}$, then \eqref{A4} has a positive radial decreasing ground state solution;
\end{itemize}
\item [$\left(2\right)$] if $p=q=2_{1}^*$, then \eqref{A4} has no nontrivial solution.
\end{itemize}
\end{thm}
The sharp classification of existence and nonexistence of solutions to \eqref{A4} relies on the value of $\alpha_{1}$ and $\lambda_{1}$, which can be estimated by \cite{Cassani-Tarsi-Zhang,Du}. We next extend the results of Theorem \ref{thm6} to the fractional case where $\alpha_{s}$ and $\lambda_{s}$ have explicit bounds by Theorem \ref{thm2}.
\begin{thm}\label{thm7}
Let $0<s<1$, $N>2s$ and $0<\lambda<1$.
\begin{itemize}
\item [$\left(1\right)$] if $2<p,q<2_{s}^*$, then \eqref{A4} has a positive ground state solution.
\item [$\left(2\right)$] if $2<p<2_{s}^*$ and $q=2_{s}^*$, let
$$
\alpha_{s}=\left[\frac{\mathcal{S}_{s}^{\frac{N}{2s}}}{\frac{N}{s}\left(\frac{1}{2}-\frac{1}{q}\right)
S_{s,q}^{\frac{q}{q-2}}\left(\mathbb{R}^{N}\right)}\right]^{\frac{1}{\frac{q}{q-2}-\frac{N}{2s}}},
$$
then there exists $\lambda_{s}\in\left[\sqrt{1-\alpha_{s}},1\right)$ such that
\begin{itemize}
\item [$\left(i\right)$] if $\lambda<\lambda_{s}$, then \eqref{A4} has no ground state solution;
\item [$\left(ii\right)$] if $\lambda>\lambda_{s}$, then \eqref{A4} has a positive radial decreasing ground state solution.
\end{itemize}
\item [$\left(3\right)$] if $p=q=2_{s}^*$, then \eqref{A4} has no nontrivial solution.
\end{itemize}
\end{thm}

\section{Preliminaries}
\subsection{Fractional Sobolev spaces and the fractional Laplacian}
Let $0<s<1\leq p<+\infty$, the so-called Gagliardo semi-norm is given by
$$
\left[u\right]_{W^{s,p}\left(\mathbb{R}^{N}\right)}
=\left(\int_{\mathbb{R}^{N}}\int_{\mathbb{R}^{N}}
\frac{\left|u\left(x\right)-u\left(y\right)\right|^{p}}{\left|x-y\right|^{N+sp}}dxdy\right)^{\frac{1}{p}}.
$$
The fractional Sobolev space $W^{s,p}\left(\mathbb{R}^{N}\right)$ is defined as the completion of $C_{0}^{\infty}\left(\mathbb{R}^{N}\right)$ with respect to the norm
$$
\left\|u\right\|_{W^{s,p}\left(\mathbb{R}^{N}\right)}=\left(\left[u\right]^{p}_{W^{s,p}\left(\mathbb{R}^{N}\right)}
+\left\|u\right\|^{p}_{L^{p}\left(\mathbb{R}^{N}\right)}\right)^{\frac{1}{p}},
$$
then $W^{s,p}_{0}\left(\Omega\right)$ is defined by
$$
W^{s,p}_{0}\left(\Omega\right):=\left\{u\in W^{s,p}\left(\mathbb{R}^{N}\right),u\equiv0~\text{in}~\mathbb{R}^{N}\backslash \Omega\right\},
$$
and $D^{s,p}\left(\Omega\right)$ is the completion of $C_{0}^{\infty}\left(\Omega\right)$ with respect to $\left[u\right]_{W^{s,p}\left(\mathbb{R}^{N}\right)}$.

In the Hilbert case $p=2$, let $\mathcal{F}$ be the standard Fourier transform
$$
\mathcal{F}u\left(\xi\right)=\int_{\mathbb{R}^{N}}u\left(x\right)e^{-2\pi i\xi\cdot x}dx,
$$
then the fractional Laplace operator $\left(-\Delta\right)^{s}$ is defined by
\begin{align}\label{D}
\left(-\Delta\right)^{s}u=\mathcal{F}^{-1}\left(\left|2\pi\xi\right|^{2s}\mathcal{F}u\right).
\end{align}
The Hilbert space $W^{s,2}\left(\mathbb{R}^{N}\right)$ coincides with the Bessel potential space $H^{s}\left(\mathbb{R}^{N}\right)$ defined via Fourier transform
$$
H^{s}\left(\mathbb{R}^{N}\right):=\left\{u\in L^{2}\left(\mathbb{R}^{N}\right),\mathcal{F}^{-1}\left[\left(1+\left|2\pi\xi\right|^{2s}\right)
^{\frac{1}{2}}\mathcal{F}u\right]\in L^{2}\left(\mathbb{R}^{N}\right)\right\}
$$
with the inner product
$$
\left<u,v\right>_{H^{s}\left(\mathbb{R}^{N}\right)}
=\int_{\mathbb{R}^{N}}\left(-\Delta\right)^{\frac{s}{2}}u\left(-\Delta\right)^{\frac{s}{2}}v+uvdx,
$$
and endowed with the norm
$$
\left\|u\right\|_{H^{s}\left(\mathbb{R}^{N}\right)}=\left<u,u\right>_{H^{s}\left(\mathbb{R}^{N}\right)}\ .
$$
Similarly $H_{0}^{s}\left(\Omega\right)$ is defined by
$$
H_{0}^{s}\left(\Omega\right):=\left\{u\in H^{s}\left(\mathbb{R}^{N}\right),u\equiv0~\text{in}~\mathbb{R}^{N}\backslash \Omega\right\}.
$$
Let us mention that for all $u\in H_{0}^{1}\left(\Omega\right)$, on the one hand by the formula of Bourgain-Br\'{e}zis-Mironescu \cite{Bourgain-Brezis-Mironescu} we have
$$
\lim_{s\rightarrow1^{-}}\left(1-s\right)\left[u\right]^{2}_{W^{s,2}\left(\mathbb{R}^{N}\right)}
=\frac{\pi^{\frac{N}{2}}}{2\Gamma\left(\frac{N}{2}+1\right)}\left\| \nabla u\right\|^{2}_{L^{2}\left(\Omega\right)} \ .
$$
On the other hand, by the formula of Maz'ya-Shaposhnikova \cite{Mazya-Shaposhnikova} we have
$$
\lim_{s\rightarrow0^{+}}s\left[u\right]^{2}_{W^{s,2}\left(\mathbb{R}^{N}\right)}
=\frac{N\pi^{\frac{N}{2}}}{\Gamma\left(\frac{N}{2}+1\right)}\left\| u\right\|^{2}_{L^{2}\left(\Omega\right)},
$$
and hence by the identity \eqref{B}, one has
\begin{align*}
\left\{
\begin{array}{lcl}
\displaystyle \lim_{s\rightarrow1^{-}}\left\|\left(-\Delta\right)^{\frac{s}{2}}u\right\|_{L^{2}\left(\mathbb{R}^{N}\right)}
=\left\| \nabla u\right\|_{L^{2}\left(\Omega\right)},\\
\\
\displaystyle \lim_{s\rightarrow0^{+}}\left\|\left(-\Delta\right)^{\frac{s}{2}}u\right\|_{L^{2}\left(\mathbb{R}^{N}\right)}
=\left\| u\right\|_{L^{2}\left(\Omega\right)},
\end{array}
\right.
\end{align*}
and thus the Fourier characterization of $H_{0}^{s}\left(\Omega\right)$ recovers both the norm of $W_{0}^{1,2}\left(\Omega\right)$ and $L^{2}\left(\Omega\right)$.

\subsection{Localizing issues}

Let us briefly discuss an alternative definition of $\mathcal{S}_{s,2_{s}^{\ast}}$ which involves the localized Gagliardo semi-norm
$$
\left[u\right]_{W^{s,2}\left(\Omega\right)}
=\left(\int_{\Omega}\int_{\Omega}\frac{\left|u\left(x\right)-u\left(y\right)\right|^{2}}
{\left|x-y\right|^{N+2s}}dxdy\right)^{\frac{1}{2}}.
$$
Let $0<s<1\leq p<+\infty$, the fractional Sobolev space $\widetilde{W}^{s,2}\left(\Omega\right)$ is defined as the completion of $C_{0}^{\infty}\left(\mathbb{R}^{N}\right)$ with respect to the norm
$$
\left\|u\right\|_{\widetilde{W}^{s,2}\left(\Omega\right)}
=\left(\left[u\right]^{2}_{W^{s,2}\left(\Omega\right)}+\left\|u\right\|^{2}_{L^{2}\left(\Omega\right)}\right)^{\frac{1}{2}},
$$
then $\widetilde{W}^{s,2}_{0}\left(\Omega\right)$ is defined by
$$
\widetilde{W}^{s,2}_{0}\left(\Omega\right):=\left\{u\in \widetilde{W}^{s,2}\left(\Omega\right),u\equiv0~in~\mathbb{R}^{N}\backslash \Omega\right\},
$$
and the space $\widetilde{D}^{s,2}\left(\Omega\right)$ is the completion of $C_{0}^{\infty}\left(\mathbb{R}^{N}\right)$ with respect to the norm $\left[u\right]_{W^{s,2}\left(\Omega\right)}$. Obviously, we have
$$
\widetilde{W}^{s,2}\left(\mathbb{R}^{N}\right)=W^{s,2}\left(\mathbb{R}^{N}\right).
$$
Moreover, Brasco-Lindgren-Parini \cite{Brasco-Lindgren-Parini} showed that $\widetilde{W}_{0}^{s,2}\left(\Omega\right)$ and $W_{0}^{s,2}\left(\Omega\right)$ do coincide if $s\neq\frac{1}{2}$.

Actually, for $N>2s$, let us consider the possible embedding
\begin{align}\label{E}
\widetilde{D}^{s,2}\left(\Omega\right)\hookrightarrow L^{2_{s}^{\ast}}\left(\Omega\right)
\end{align}
and denote $\widetilde{S}_{s,2_{s}^{\ast}}\left(\Omega\right)$ as the optimal constant such that
$$
\widetilde{S}_{s,2_{s}^{\ast}}\left(\Omega\right)\left\|u\right\|^{2}_{L^{2^{\ast}_{s}}\left(\Omega\right)}
\leq\left[u\right]_{W^{s,2}\left(\Omega\right)}^{2},\quad u\in \widetilde{D}^{s,2}\left(\Omega\right)\backslash\left\{0\right\},
$$
namely
$$
\widetilde{S}_{s,2_{s}^{\ast}}\left(\Omega\right)
=\inf_{u\in \widetilde{D}^{s,2}\left(\Omega\right)\backslash\left\{0\right\}}
\frac{\left[u\right]^{2}_{W^{s,2}\left(\Omega\right)}}{\left\|u\right\|^{2}_{L^{2^{\ast}_{s}}\left(\Omega\right)}}.
$$
Nevertheless, in contrast to $S_{s,2_{s}^{\ast}}\left(\Omega\right)$, there is no scale invariance for $\widetilde{S}_{s,2_{s}^{\ast}}\left(\Omega\right)$, which means that $\widetilde{S}_{s,2_{s}^{\ast}}\left(\Omega\right)$ strictly depends on $\Omega$. Moreover, when $0<s<\frac{1}{2}$, Frank-Jin-Xiong \cite{Frank-Jin-Xiong} proved that $\widetilde{S}_{s,2_{s}^{\ast}}\left(\Omega\right)=0$, which implies \eqref{E} fails. When $N\geq2$ and $\frac{1}{2}<s<1$, the constant $\widetilde{S}_{s,2_{s}^{\ast}}\left(\Omega\right)$ can be achieved provided some additional conditions are assumed as done in Frank-Jin-Xiong \cite{Frank-Jin-Xiong}. Furthermore, Dyda-Frank \cite{Dyda-Frank} showed that there exists a uniform constant
$$
\widetilde{\mathcal{S}}_{s}=\inf_{\Omega\subset\mathbb{R}^{N}}\widetilde{S}_{s,2_{s}^{\ast}}\left(\Omega\right),
$$
such that for any $\Omega\neq\mathbb{R}^{N}$, there holds
$$
\widetilde{\mathcal{S}}_{s}\left\|u\right\|^{2}_{L^{2^{\ast}_{s}}\left(\Omega\right)}
\leq\left[u\right]_{W^{s,2}\left(\Omega\right)}^{2},
\quad u\in \widetilde{D}^{s,2}\left(\Omega\right)\backslash\left\{0\right\}.
$$
One possible explanation for such phenomena goes back to the Br\'{e}zis-Nirenberg \cite{Brezis-Nirenberg} result. Indeed, let us rewrite the norm
$$
\left[u\right]_{W^{s,2}\left(\Omega\right)}=\left[u\right]_{W^{s,2}\left(\mathbb{R}^{N}\right)}
-2\int_{\mathbb{R}^{N}-\Omega}\int_{\Omega}\frac{\left|u\left(x\right)\right|^{2}}{\left|x-y\right|^{N+2s}}dxdy,
$$
so that we see how the negative part in the right hand side lowers the value $\mathcal{S}_{s,2_{s}^{\ast}}$ and as a consequence $\widetilde{S}_{s,2_{s}^{\ast}}\left(\Omega\right)$ retrives a minimizer.

\section{Bounds for best constants of fractional subcritical Sobolev embeddings}

In this section, we establish fine bounds for $S_{s,q}\left(\Omega\right)$ and $S_{s,q}\left(\mathbb{R}^{N}\right)$ in the borderline case $p=1$, the Hilbert case $2=p<\frac{N}{s}$ and the limiting case $s=\frac{1}{2}$, $p=2$ and $N=1$. Moreover, we also establish sharp asymptotics for the limiting case $s=\frac{1}{2}$, $p=2$ and $N=1$.

\subsection{The borderline case $p=1$: proof of Theorem \ref{thm1}}

We first prove Theorem \ref{thm1}-$\left(1\right)$. By H\"{o}lder's inequality, we have
\begin{align*}
S_{s,q}\left(\Omega\right)
\geq\left|\Omega\right|^{\frac{1}{1^{\ast}_{s}}-\frac{1}{q}}
\inf_{u\in W^{s,1}_{0}\left(\Omega\right)\backslash\left\{0\right\}}
\frac{\left[u\right]_{W^{s,1}\left(\mathbb{R}^{N}\right)}}{\left\|u\right\|_{L^{1^{\ast}_{s}}\left(\Omega\right)}}
=\mathcal{S}_{s,1_{s}^{\ast}}\left|\Omega\right|^{\frac{1}{1^{\ast}_{s}}-\frac{1}{q}}.
\end{align*}
Let us take the characteristic function $\chi_{B_{1}}$ to get
\begin{align*}
S_{s,q}\left(B_{1}\left(0\right)\right)
\leq\omega_{N}^{\frac{1}{1^{\ast}_{s}}-\frac{1}{q}}
\frac{\left[\chi_{B_{1}}\right]_{W^{s,1}\left(\mathbb{R}^{N}\right)}}{\left\|\chi_{B_{1}}\right\|_{L^{1^{\ast}_{s}}
\left(B_{1}\left(0\right)\right)}}
=\omega_{N}^{\frac{1}{1^{\ast}_{s}}-\frac{1}{q}} \mathcal{S}_{s,1_{s}^{\ast}}.
\end{align*}
Next translate the center of $B_{R_{\Omega}}$ into the origin and we apply the dilation group action \eqref{A} to get
\begin{align*}
S_{s,q}\left(\Omega\right)\leq\omega_{N}^{\frac{1}{1^{\ast}_{s}}-\frac{1}{q}} \mathcal{S}_{s,1_{s}^{\ast}}
R_{\Omega}^{N\left(\frac{1}{1^{\ast}_{s}} -\frac{1}{q}\right)}=\mathcal{S}_{s,1_{s}^{\ast}}\left|B_{R_{\Omega}}\right|^{\frac{1}{1^{\ast}_{s}}-\frac{1}{q}}.
\end{align*}
Next we prove Theorem \ref{thm1}-$\left(3\right)$. For any $u\in W^{s,1}\left(\mathbb{R}^{N}\right)\backslash\left\{0\right\}$, by interpolation inequality, we have
\begin{align*}
\left\|u\right\|_{L^{q}\left(\mathbb{R}^{N}\right)}
=\left\|u\right\|^{\lambda_{1}}_{L^{1}\left(\mathbb{R}^{N}\right)}
\left\|u\right\|^{1-\lambda_{1}}_{L^{1^{\ast}_{s}}\left(\mathbb{R}^{N}\right)},
\end{align*}
where
$$
\lambda_{1}=\left(\frac{N}{s}-1\right)\left(\frac{1^{\ast}_{s}}{q}-1\right).
$$
By Young's inequality, we get
\begin{align*}
\left\|u\right\|^{\lambda_{1}}_{L^{1}\left(\mathbb{R}^{N}\right)}\left\|u\right\|^{1-\lambda_{1}}_{L^{1^{\ast}_{s}}\left(\mathbb{R}^{N}\right)}
&=\left(\varepsilon_{1}\left\|u\right\|^{\lambda_{1}}_{L^{1}\left(\mathbb{R}^{N}\right)}\right)
\left(\frac{1}{\varepsilon_{1}}\left\|u\right\|^{1-\lambda_{1}}_{L^{1^{\ast}_{s}}\left(\mathbb{R}^{N}\right)}\right)\\
&\leq\lambda_{1}\varepsilon_{1}^{\frac{1}{\lambda_{1}}}\left\|u\right\|_{L^{1}\left(\mathbb{R}^{N}\right)}
+\left(1-\lambda_{1}\right)\varepsilon_{1}^{\frac{1}{\lambda_{1}-1}}\left\|u\right\|_{L^{1^{\ast}_{s}}\left(\mathbb{R}^{N}\right)}.
\end{align*}
Let us choose
\begin{align*}
\left\{
\begin{array}{lcl}
\displaystyle \lambda_{1}\varepsilon_{1}^{\frac{1}{\lambda_{1}}}=\rho_{1}; \\
\\
\displaystyle \left(1-\lambda_{1}\right)\varepsilon_{1}^{\frac{1}{\lambda_{1}-1}}=\rho_{1}\mathcal{S}_{s,1_{s}^{\ast}},
\end{array}
\right.
\end{align*}
hence
$$
\left\|u\right\|_{L^{q}\left(\mathbb{R}^{N}\right)}
\leq\rho_{1}\left(\left\|u\right\|_{L^{1}\left(\mathbb{R}^{N}\right)}
+\mathcal{S}_{s,1_{s}^{\ast}}\left\|u\right\|_{L^{1^{\ast}_{s}}\left(\mathbb{R}^{N}\right)}\right),
$$
where
\begin{align*}
\left\{
\begin{array}{lcl}
\displaystyle \varepsilon_{1}=\left(\frac{\lambda_{1}\mathcal{S}_{s,1_{s}^{\ast}}}{1-\lambda_{1}}\right)^{\lambda_{1}\left(\lambda_{1}-1\right)}; \\
\\
\displaystyle \rho_{1}=\lambda_{1}^{\lambda_{1}}\left(\frac{\mathcal{S}_{s,1_{s}^{\ast}}}{1-\lambda_{1}}\right)^{\lambda_{1}-1}.
\end{array}
\right.
\end{align*}
Therefore we have
\begin{align*}
\left\|u\right\|_{W^{s,1}\left(\mathbb{R}^{N}\right)}
\geq\frac{1}{\rho_{1}}\left\|u\right\|_{L^{q}\left(\mathbb{R}^{N}\right)},
\end{align*}
and we conclude
\begin{align*}
S_{s,q}\left(\mathbb{R}^{N}\right)
\geq\left[\left(\frac{N}{s}-1\right)\left(\frac{1^{\ast}_{s}}{q}-1\right)\right]^{\left(\frac{N}{s}-1\right)
\left(1-\frac{1^{\ast}_{1}}{q}\right)}\left[\frac{N}{s}\left(1-\frac{1}{q}\right)\right]
^{\frac{N}{s}\left(\frac{1}{q}-1\right)}
\mathcal{S}_{s,1_{s}^{\ast}}^{\frac{N}{s}\left(1-\frac{1}{q}\right)}.
\end{align*}
Let us choose the characteristic function $\chi_{B_{k}}$ to have
\begin{align*}
S_{s,q}\left(\mathbb{R}^{N}\right)
&\leq\left(\omega_{N}k^{N}\right)^{\frac{1}{1_{s}^{\ast}}-\frac{1}{q}}\frac{\left[\chi_{B_{k}}\right]_{W^{s,1}
\left(\mathbb{R}^{N}\right)}}{\left\|\chi_{B_{k}}\right\|_{L^{1_{s}^{\ast}}\left(\mathbb{R}^{N}\right)}}
+\frac{\left\|\chi_{B_{k}}\right\|_{L^{1}
\left(\mathbb{R}^{N}\right)}}{\left\|\chi_{B_{k}}\right\|_{L^{q}\left(\mathbb{R}^{N}\right)}}\\
&=\omega_{N}^{\frac{1}{1_{s}^{\ast}}-\frac{1}{q}}\mathcal{S}_{s,1_{s}^{\ast}}k^{N\left(\frac{1}{1_{s}^{\ast}}-\frac{1}{q}\right)}
+\omega_{N}^{1-\frac{1}{q}}k^{N\left(1-\frac{1}{q}\right)}\\
&=:g_{1}\left(k\right).
\end{align*}
We conclude the following
\begin{align*}
S_{1,q}\left(\mathbb{R}^{N}\right)
&\leq \inf_{k>0}g_{1}\left(k\right)\\
&=g_{1}\left\{\left[\frac{\mathcal{S}_{s,1_{s}^{\ast}}\left(\frac{1}{q}-\frac{1}{1_{s}^{\ast}}\right)}
{\omega_{N}^{\frac{s}{N}}\left(1-\frac{1}{q}\right)}\right]^{\frac{1}{s}}\right\}\\
&=\frac{s}{N}\left(\frac{1}{q}-\frac{1}{1_{s}^{\ast}}\right)^{\frac{N}{s}\left(\frac{1}{1_{s}^{\ast}}-\frac{1}{q}\right)}
\left(1-\frac{1}{q}\right)^{\frac{N}{s}\left(\frac{1}{q}-1\right)}\mathcal{S}_{s,1_{s}^{\ast}}^{\frac{N}{s}\left(1-\frac{1}{q}\right)}.
\end{align*}
Finally, when $q=1$, there holds
\begin{align*}
S_{s,1}\left(\mathbb{R}^{N}\right)
\leq\lim_{k\rightarrow+\infty}g_{1}\left(k\right)
=1,
\end{align*}
therefore $S_{s,1}\left(\mathbb{R}^{N}\right)=1$, which yields Theorem \ref{thm1}-$\left(2\right)$.

\subsection{The Hilbert case $2=p<\frac{N}{s}$: proof of Theorem \ref{thm2}}

We start by proving Theorem \ref{thm2}-$\left(1\right)$. By H\"{o}lder's inequality we have
\begin{align*}
S_{s,q}\left(\Omega\right)
\geq\left|\Omega\right|^{2\left(\frac{1}{2^{\ast}_{s}}-\frac{1}{q}\right)}\inf_{u\in H^{s}_{0}\left(\Omega\right)\backslash\left\{0\right\}}
\frac{\left\|\left(-\Delta\right)^{\frac{s}{2}}u\right\|^{2}_{L^{2}\left(\mathbb{R}^{N}\right)}}{\left\|u\right\|^{2}_{L^{2^{\ast}_{s}}
\left(\Omega\right)}}
=\mathcal{S}_{s}\left|\Omega\right|^{2\left(\frac{1}{2^{\ast}_{s}}-\frac{1}{q}\right)}.
\end{align*}
Next we establish an upper bound for $S_{1,q}\left(\Omega\right)$ by using the classical Fourier transform of radial functions by Stein-Weiss \cite{Stein-Weiss} involving the standard Bessel function
\begin{align*}
\mathrm{J}_{v}\left(t\right)=\left\{
\begin{array}{lcl}
\displaystyle \frac{\left(\frac{t}{2}\right)^{v}}{\Gamma\left(v+\frac{1}{2}\right)\Gamma\left(\frac{1}{2}\right)}
\int^{1}_{-1}\left(1-s^{2}\right)^{v-\frac{1}{2}}e^{its}ds,&& v>-\frac{1}{2}; \\
\\
\displaystyle \sqrt{\frac{2}{\pi}}\frac{\cos t}{\sqrt{t}}, && v=-\frac{1}{2}.
\end{array}
\right.
\end{align*}
Consider the radial decreasing function
\begin{align*}
f_{1}\left(x\right)=\left\{
\begin{array}{lcl}
\left(k^{2}-\left|x\right|^{2}\right)^{s}, && 0\leq \left|x\right|\leq k, \\
\\
0,&& k\leq\left|x\right|\leq 1,
\end{array}
\right.
\end{align*}
we have
\begin{align*}
\mathcal{F}f_{1}\left(\xi\right)
&=2\pi\left|\xi\right|^{-\frac{N}{2}+1}
\int^{+\infty}_{0}f_{1}\left(r\right)r^{\frac{N}{2}}\mathrm{J}_{\frac{N}{2}-1}\left(2\pi \left|\xi\right|r\right)dr\\
&=\pi^{-s}\left|\xi\right|^{-\frac{N}{2}-s}k^{\frac{N}{2}+s}\Gamma\left(s+1\right)\mathrm{J}_{\frac{N}{2}+s}\left(2\pi k \left|\xi\right|\right),
\end{align*}
where we apply the following identity
$$
\int^{1}_{0}\left(1-r^{2}\right)^{s}r^{\frac{N}{2}}\mathrm{J}_{\frac{N}{2}-1}\left(2\pi \left|\xi\right|kr\right)dr=2^{-1}\left(\pi\left|\xi\right|k\right)^{-s-1}\Gamma\left(s+1\right)\mathrm{J}_{\frac{N}{2}+s}\left(2\pi \left|\xi\right|k\right).
$$
Hence \eqref{D} implies
\begin{align*}
\left\|\left(-\Delta\right)^{\frac{s}{2}}f_{1}\right\|^{2}_{L^{2}\left(\mathbb{R}^{N}\right)}
=\frac{2^{2s}\omega_{N}N}{N+2s}\left[\Gamma\left(s+1\right)\right]^{2}k^{N+2s}.
\end{align*}
Moreover, from
\begin{align*}
\left\|f_{1}\right\|^{q}_{L^{q}\left(B_{1}\left(0\right)\right)}
=\frac{\omega_{N}N}{2}\mathrm{B}\left(\frac{N}{2},qs+1\right)k^{N+2qs},
\end{align*}
we conclude that
\begin{align*}
S_{s,q}\left(B_{1}\left(0\right)\right)
\leq \frac{2^{2s+\frac{2}{q}}\left(\omega_{N}N\right)^{1-\frac{2}{q}}}{N+2s}
\left[\Gamma\left(s+1\right)\right]^{2}\left[\mathrm{B}
\left(\frac{N}{2},qs+1\right)\right]^{-\frac{2}{q}}k^{2N\left(\frac{1}{2^{\ast}_{s}}-\frac{1}{q}\right)}.
\end{align*}
We reach the desired result letting $k\rightarrow1$. By using \eqref{A} again we get
\begin{align*}
S_{s,q}\left(\Omega\right)\leq\frac{2^{2s+\frac{2}{q}}\left(\omega_{N}N\right)^{1-\frac{2}{q}}}{N+2s}
\left[\Gamma\left(s+1\right)\right]^{2}\left[\mathrm{B}
\left(\frac{N}{2},qs+1\right)\right]^{-\frac{2}{q}}
R_{\Omega}^{2N\left(\frac{1}{2^{\ast}_{s}} -\frac{1}{q}\right)}.
\end{align*}
Next we prove Theorem \ref{thm2}-$\left(3\right)$. For any $u\in H^{s}\left(\mathbb{R}^{N}\right)\backslash\left\{0\right\}$, by interpolation inequality, we have
\begin{align*}
\left\|u\right\|^{2}_{L^{q}\left(\mathbb{R}^{N}\right)}
=\left\|u\right\|^{2\lambda_{2}}_{L^{2}\left(\mathbb{R}^{N}\right)}
\left\|u\right\|^{2\left(1-\lambda_{2}\right)}_{L^{2^{\ast}_{s}}\left(\mathbb{R}^{N}\right)},
\end{align*}
where
$$
\lambda_{2}=\frac{N}{s}\left(\frac{1}{q}-\frac{1}{2_{s}^{\ast}}\right).
$$
By Young's inequality, we get
\begin{align*}
\left\|u\right\|^{2\lambda_{2}}_{L^{2}\left(\mathbb{R}^{N}\right)}
\left\|u\right\|^{2\left(1-\lambda_{2}\right)}_{L^{2^{\ast}_{s}}\left(\mathbb{R}^{N}\right)}
&=\left(\varepsilon_{2}\left\|u\right\|^{2\lambda_{2}}_{L^{2}\left(\mathbb{R}^{N}\right)}\right)
\left(\frac{1}{\varepsilon_{2}}\left\|u\right\|^{2\left(1-\lambda_{2}\right)}_{L^{2^{\ast}_{s}}\left(\mathbb{R}^{N}\right)}\right)\\
&\leq\lambda_{2}\varepsilon_{2}^{\frac{1}{\lambda_{2}}}\left\|u\right\|^{2}_{L^{2}\left(\mathbb{R}^{N}\right)}
+\left(1-\lambda_{2}\right)\varepsilon_{2}^{\frac{1}{\lambda_{2}-1}}\left\|u\right\|^{2}_{L^{2^{\ast}_{s}}\left(\mathbb{R}^{N}\right)}.
\end{align*}
Let us set
\begin{align*}
\left\{
\begin{array}{lcl}
\displaystyle \lambda_{2}\varepsilon_{2}^{\frac{1}{\lambda_{2}}}=\rho_{2}; \\
\\
\displaystyle \left(1-\lambda_{2}\right)\varepsilon_{2}^{\frac{1}{\lambda_{2}-1}}=\rho_{2}\mathcal{S}_{s},
\end{array}
\right.
\end{align*}
hence
$$
\left\|u\right\|^{2}_{L^{q}\left(\mathbb{R}^{N}\right)}
\leq\rho_{2}\left(\left\|u\right\|^{2}_{L^{2}\left(\mathbb{R}^{N}\right)}
+\mathcal{S}_{s}\left\|u\right\|^{2}_{L^{2^{\ast}_{s}}\left(\mathbb{R}^{N}\right)}\right),
$$
where
\begin{align*}
\left\{
\begin{array}{lcl}
\displaystyle \varepsilon_{2}=\left(\frac{\lambda_{2}\mathcal{S}_{s}}{1-\lambda_{2}}\right)^{\lambda_{2}\left(\lambda_{2}-1\right)}; \\
\\
\displaystyle \rho_{2}=\lambda_{2}^{\lambda_{2}}\left(\frac{\mathcal{S}_{s}}{1-\lambda_{2}}\right)^{\lambda_{2}-1}.
\end{array}
\right.
\end{align*}
Therefore one has
\begin{align*}
\left\|u\right\|^{2}_{H^{s}\left(\mathbb{R}^{N}\right)}
\geq\frac{1}{\rho_{2}}\left\|u\right\|^{2}_{L^{q}\left(\mathbb{R}^{N}\right)},
\end{align*}
from which we deduce
\begin{align*}
S_{s,q}\left(\mathbb{R}^{N}\right)
\geq\left[\frac{N}{s}\left(\frac{1}{q}-\frac{1}{2_{s}^{\ast}}\right)\right]^{\frac{N}{s}\left(\frac{1}{2_{s}^{\ast}}-\frac{1}{q}\right)}
\left[\frac{N}{s}\left(\frac{1}{2}-\frac{1}{q}\right)\right]
^{\frac{N}{s}\left(\frac{1}{q}-\frac{1}{2}\right)}
\mathcal{S}_{s}^{\frac{N}{s}\left(\frac{1}{2}-\frac{1}{q}\right)}.
\end{align*}
Let us now consider the following test function
\begin{align*}
f_{2}\left(x\right)=\left\{
\begin{array}{lcl}
\left(k^{2}-\left|x\right|^{2}\right)^{s}, && 0\leq \left|x\right|< k, \\
\\
0, && \left|x\right|\geq k.
\end{array}
\right.
\end{align*}
and direct calculations as for $f_1$ give
\begin{align*}
S_{s,q}\left(\mathbb{R}^{N}\right)
\leq&\left(\omega_{N}N\right)^{1-\frac{2}{q}}\left\{\frac{2^{2s+\frac{2}{q}}}{N+2s}
\Gamma^{2}\left(s+1\right)
\left[\mathrm{B}\left(\frac{N}{2},qs+1\right)\right]^{-\frac{2}{q}}k^{2N\left(\frac{1}{2^{\ast}_{s}}-\frac{1}{q}\right)}\right.\\
&+\left.2^{\frac{2}{q}-1}\mathrm{B}
\left(\frac{N}{2},2s+1\right)\left[\mathrm{B}\left(\frac{N}{2},qs+1\right)\right]
^{-\frac{2}{q}}k^{N\left(1-\frac{2}{q}\right)}\right\}\\
=:&g_{2}\left(k\right).
\end{align*}
We conclude that
\begin{align*}
S_{s,q}\left(\mathbb{R}^{N}\right)
\leq&\inf_{k>0}g_{2}\left(k\right)\\
=&g_{2}\left\{\left[\frac{2^{2s+1}
\Gamma^{2}\left(s+1\right)\left(\frac{1}{q}-\frac{1}{2^{\ast}_{s}}\right)}{\left(N+2s\right)\mathrm{B}\left(\frac{N}{2},2s+1\right)
\left(\frac{1}{2}-\frac{1}{q}\right)}\right]^{\frac{1}{2s}}\right\}\\
=& 2^{\left(2N+\frac{N}{s}-2\right)\left(\frac{1}{2}-\frac{1}{q}\right)}\omega^{1-\frac{2}{q}}_{N}s
\left[N\mathrm{B}\left(\frac{N}{2},qs+1\right)\right]^{-\frac{2}{q}}\\
&\cdot\left[\frac{\Gamma^{2}\left(s+1\right)}{\left(N+2s\right)
\left(\frac{1}{2}-\frac{1}{q}\right)}\right]^{\frac{N}{s}\left(\frac{1}{2}-\frac{1}{q}\right)}
\left[\frac{\mathrm{B}\left(\frac{N}{2},2s+1\right)}{\frac{1}{q}-\frac{1}{2^{\ast}_{s}}}\right]
^{\frac{N}{s}\left(\frac{1}{q}-\frac{1}{2^{\ast}_{s}}\right)}.
\end{align*}
Finally, when $q=2$ we have
$$
S_{s,2}\left(\mathbb{R}^{N}\right)\leq\lim_{k\rightarrow+\infty}g_{2}\left(k\right)=1,
$$
therefore $S_{s,2}\left(\mathbb{R}^{N}\right)=1$, which yields Theorem \ref{thm2}-$\left(2\right)$.

\subsection{The limiting case $s=\frac{1}{2}$, $p=2$ and $N=1$: proof of Theorem \ref{thm3}}
Let us prove Theorem \ref{thm3}-$\left(1\right)$. Consider the following so-called Moser-type function defined in the interval $\left[-1,1\right]$ by
\begin{align*}
f_{3}\left(x\right)=\left\{
\begin{array}{lcl}
\ln K-\ln k, && 0\leq\left|x\right|\leq k; \\
\\
\ln K-\ln\left|x\right|,&& k\leq\left|x\right|\leq K;\\
\\
0,&&K\leq\left|x\right|\leq 1.
\end{array}
\right.
\end{align*}
A direct computation gives
\begin{align*}
\mathcal{F}f_{3}\left(\xi\right)
=\frac{1}{\pi\xi}\int^{K}_{ k}\frac{\sin \left(2\pi\xi x\right)}{x}dx.
\end{align*}
Hence, by \eqref{D} and Fubini's theorem, we get
\begin{align*}
\left\|\left(-\Delta\right)^{\frac{1}{4}}f_{3}\right\|^{2}_{L^{2}\left(\mathbb{R}\right)}
&=\frac{4}{\pi}\int^{K}_{ k}\int^{K}_{ k}\frac{1}{xy}\left[\int_{0}^{+\infty}\frac{\sin \left(2\pi x\xi\right)\sin \left(2\pi y\xi\right)}{\xi}d\xi\right]dxdy\\
&=\frac{2}{\pi}\int^{K}_{ k}\int^{K}_{ k}\frac{1}{xy}\ln\left|\frac{x+y}{x-y}\right|dxdy\\
&=\frac{2}{\pi}\int^{K}_{ k}\frac{1}{y}\left(\int^{\frac{K}{y}}_{ \frac{k}{y}}\frac{1}{t}\ln\left|\frac{t+1}{t-1}\right|dt\right)dy.
\end{align*}
Notice that
\begin{align*}
\int^{\frac{K}{y}}_{ \frac{k}{y}}\frac{1}{t}\ln\left|\frac{t+1}{t-1}\right|dt\leq\int^{+\infty}_{ 0}\frac{1}{t}\ln\left|\frac{t+1}{t-1}\right|dt=\frac{\pi^{2}}{2}
\end{align*}
and thus
\begin{align*}
\left\|\left(-\Delta\right)^{\frac{1}{4}}f_{3}\right\|^{2}_{L^{2}\left(\mathbb{R}\right)}\leq\pi\left(\ln K-\ln k\right).
\end{align*}
Moreover, we have
\begin{align*}
\left\|f_{3}\right\|^{q}_{L^{q}\left(\left[-1,1\right]\right)}
\geq2k\left(\ln K-\ln k\right)^{q},
\end{align*}
and hence there holds
\begin{align*}
S_{\frac{1}{2},q}\left(\left[-1,1\right]\right)
&\leq2^{-\frac{2}{q}}\pi k^{-\frac{2}{q}}\left(\ln K-\ln k\right)^{-1}\\
&=:g_{3}\left(k,K\right).
\end{align*}
We know that
\begin{align*}
S_{\frac{1}{2},q}\left(B_{1}\left(0\right)\right)
&\leq \inf_{0<k<K\leq1}g_{3}\left(k,K\right)\\
&=g_{3}\left(e^{-\frac{q}{2}},1\right)\\
&=\frac{2^{1-\frac{2}{q}}\pi e}{q}.
\end{align*}
By \eqref{A}, we have
\begin{align*}
S_{\frac{1}{2},q}\left(\Omega\right)\leq\frac{2^{1-\frac{2}{q}}\pi e}{q}
R_{\Omega}^{-\frac{2}{q}}.
\end{align*}
Next we prove Theorem \ref{thm3}-$\left(3\right)$. Let us consider
\begin{align*}
f_{4}\left(x\right)=\left\{
\begin{array}{lcl}
\ln K-\ln k, && 0\leq\left|x\right|\leq k; \\
\\
\ln K-\ln\left|x\right|,&& k\leq\left|x\right|\leq K;\\
\\
0,&&\left|x\right|\geq K.
\end{array}
\right.
\end{align*}
A similar computation as for $f_{3}$, we get
\begin{align*}
S_{\frac{1}{2},q}\left(\mathbb{R}\right)
&\leq2^{-\frac{2}{q}}\left[\pi k^{-\frac{2}{q}}\left(\ln K-\ln k\right)^{-1}+2k^{-\frac{2}{q}}K\right]\\
&=:g_{4}\left(k,K\right).
\end{align*}
We conclude that
\begin{align*}
S_{\frac{1}{2},q}\left(\mathbb{R}\right)
&\leq\inf_{0<k<K}g_{4}\left(k,K\right)\\
&=g_{4}\left[\frac{2\pi}{\left(q-2\right)^{2}e^{\frac{q-2}{2}}},\frac{2\pi}{\left(q-2\right)^{2}}\right]\\
&=2^{1-\frac{4}{q}}\pi^{1-\frac{2}{q}}q\left(q-2\right)^{\frac{4}{q}-2}e^{\frac{q-2}{q}}.
\end{align*}
When $q=2$, there holds
\begin{align*}
S_{\frac{1}{2},2}\left(\mathbb{R}\right)
\leq\lim_{K\rightarrow+\infty}g_{4}\left(Ke^{-\frac{1}{\sqrt{K}}},K\right)
=1,
\end{align*}
therefore $S_{\frac{1}{2},2}\left(\mathbb{R}\right)=1$, which yields Theorem \ref{thm3}-$\left(2\right)$.

Finally, we prove Theorem \ref{thm3}-$\left(4\right)$. The bounds we established in Theorem \ref{thm3} do not obviously show the asymptotic behavior of $S_{s,q}\left(\Omega\right)$ and $S_{s,q}\left(\mathbb{R}^{N}\right)$. Actually, Ren-Wei \cite{Ren-Wei} used the Trudinger-Moser inequality
$$
\int_{\Omega}\exp\left[4\pi\left(\frac{u}{\left\|\nabla u\right\|_{L^{2}\left(\Omega\right)}}\right)^{2}\right]dx
\leq C\left|\Omega\right|
$$
to obtain \eqref{C}. Motivated by Ren-Wei \cite{Ren-Wei}, we recall the fractional Trudinger-Moser inequality on the bounded domain $\Omega$ established by Martinazzi \cite{Martinazzi}.
\begin{Prop}\label{Prop1}
Let $N=1$, for any $u\in H^{\frac{1}{2}}_{0}\left(\Omega\right)$ and $0<\gamma\leq\pi$, there exists a positive constant $C_{1}$ such that
$$
\sup_{\left\|\left(-\Delta\right)^{\frac{1}{4}}u\right\|_{L^{2}\left(\mathbb{R}\right)}\leq1}\int_{\Omega}e^{\gamma u^{2}}dx
\leq C_{1}\left|\Omega\right|.
$$
\end{Prop}
Another version of the fractional Trudinger-Moser inequality in the whole space $\mathbb{R}^{N}$ was obtained by Iula-Maalaoui-Martinazzi \cite{Iula-Maalaoui-Martinazzi}.
\begin{Prop}\label{Prop2}
For any $u\in H^{\frac{1}{2}}\left(\mathbb{R}\right)$ and $0<\gamma\leq\pi$, there exists a positive constant $C_{2}$ such that
$$
\sup_{\left\|u\right\|_{H^{\frac{1}{2}}\left(\mathbb{R}\right)}\leq1}\int_{\mathbb{R}}e^{\gamma u^{2}}-1dx\leq C_{2}.
$$
\end{Prop}

On the one hand, for any $u\in H_{0}^{\frac{1}{2}}\left(\Omega\right)$, we have
\begin{align*}
\left\|u\right\|^{q}_{L^{q}\left(\Omega\right)}
&=\pi^{-\frac{q}{2}}\left\|\left(-\Delta\right)^{\frac{1}{4}}u\right\|^{q}_{L^{2}\left(\mathbb{R}\right)}\int_{\Omega}
\left[\pi \left(\frac{u}{\left\|\left(-\Delta\right)^{\frac{1}{4}}u\right\|_{L^{2}\left(\mathbb{R}\right)}}\right)^{2}\right]^{\frac{q}{2}}dx\\
&\leq C_{1}\pi^{-\frac{q}{2}}\Gamma\left(\frac{q}{2}
+1\right)\left|\Omega\right|\left\|\left(-\Delta\right)^{\frac{1}{4}}u\right\|^{q}_{L^{2}\left(\mathbb{R}\right)},
\end{align*}
where we apply Proposition \ref{Prop1} and the inequality
$$
x^{t}\leq\Gamma\left(t+1\right)e^{x},\quad x,t\geq0.
$$
Hence we get
\begin{align*}
S_{\frac{1}{2},q}\left(\Omega\right)
\geq C_{1}^{-\frac{2}{q}}\pi\left[\Gamma\left(\frac{q}{2}+1\right)\right]^{-\frac{2}{q}}\left|\Omega\right|^{-\frac{2}{q}}.
\end{align*}
By Stirling's formula
$$
\Gamma\left(t+1\right)\sim\sqrt{2\pi}t^{t+\frac{1}{2}}e^{-t},\quad t\rightarrow+\infty
$$
and Theorem \ref{thm3}-$\left(1\right)$, we get the first asymptotic behavior
\begin{align*}
\lim_{q\rightarrow +\infty } qS_{\frac{1}{2},q}\left(\Omega\right)=2\pi e.
\end{align*}
On the other hand, let us consider the symmetric decreasing rearrangement $u^{\#}$. By the following fractional P\'{o}lya-Szeg\"{o} inequality proved by Almgren-Lieb \cite{Almgren-Lieb}:
\begin{align}\label{F}
\left[u^{\#}\right]_{W^{\frac{1}{2},2}\left(\mathbb{R}^{N}\right)}\leq\left[u\right]_{W^{\frac{1}{2},2}\left(\mathbb{R}^{N}\right)},
\end{align}
we can replace $u$ by $u^{\#}$. Let us split the norm $\left\|u^{\#}\right\|^{q}_{L^{q}\left(\mathbb{R}\right)}$ into two parts
\begin{align}\label{G}
\left\|u^{\#}\right\|^{q}_{L^{q}\left(\mathbb{R}\right)}
=\left\|u^{\#}\right\|^{q}_{L^{q}\left(\left|x\right|\leq 1\right)}
+\left\|u^{\#}\right\|^{q}_{L^{q}\left(\left|x\right|\geq 1\right)}.
\end{align}
By Proposition \ref{Prop2}, for any $u\in H^{\frac{1}{2}}\left(\mathbb{R}\right)$, we have
\begin{align}\label{H}
\left\|u^{\#}\right\|^{q}_{L^{q}\left(\left|x\right|\leq 1\right)}
&=\pi^{-\frac{q}{2}}\left\|u^{\#}\right\|^{q}_{H^{\frac{1}{2}}\left(\mathbb{R}\right)}
\int_{\left|x\right|\leq 1}\left[\pi \left(\frac{u^{\#}}{\left\|u^{\#}\right\|_{H^{\frac{1}{2}}\left(\mathbb{R}\right)}}\right)^{2}\right]^{\frac{q}{2}}dx\nonumber\\
&\leq\left(C _{2}+2\right)\pi^{-\frac{q}{2}}\Gamma\left(\frac{q}{2}+1\right)\left\|u^{\#}\right\|^{q}_{H^{\frac{1}{2}}\left(\mathbb{R}\right)}.
\end{align}
Since
$$
u^{\#}\left(r\right)\leq\frac{\left\|u^{\#}\right\|_{L^{2}\left(\mathbb{R}\right)}}{\sqrt{2r}},
$$
one has
\begin{align}\label{I}
\left\|u^{\#}\right\|^{q}_{L^{q}\left(\left|x\right|\geq 1\right)}
\leq\frac{2^{2-\frac{q}{2}}}{q-2}\left\|u^{\#}\right\|^{q}_{H^{\frac{1}{2}}\left(\mathbb{R}\right)}.
\end{align}
Let us combine \eqref{F}, \eqref{G}, \eqref{H} and \eqref{I} to obtain
\begin{align*}
S_{\frac{1}{2},q}\left(\mathbb{R}\right)
\geq\left[\left(C _{2}+2\right)\pi^{-\frac{q}{2}}\Gamma\left(\frac{q}{2}+1\right)+\frac{2^{2-\frac{q}{2}}}{q-2}\right]^{-\frac{2}{q}}.
\end{align*}
By Stirling's formula and Theorem \ref{thm3}-$\left(3\right)$, we get the second asymptotic behavior
$$
\lim_{q\rightarrow+\infty}qS_{\frac{1}{2},q}\left(\mathbb{R}\right)=2\pi e.
$$

\section{Applications to nonlocal PDEs}
\subsection{A nonlocal nonlinear Schr\"{o}dinger equation: proof of Theorem \ref{thm4}}
Let us prove Theorem \ref{thm4}-$\left(1\right)$. Define the functional $E\in \mathcal{C}^{1}\left(H^{s}\left(\mathbb{R}^{N}\right),\mathbb{R}\right)$ by
$$
E\left(u\right)=\frac{1}{2}\int_{\mathbb{R}^{N}}
\left|\left(-\Delta\right)^{\frac{s}{2}}u\right|^{2}+\left|u\right|^{2}dx
-\frac{1}{q}\int_{\mathbb{R}^{N}}Q\left|u^{+}\right|^{q}dx
$$
for which
$$
\left<E'\left(u\right),v\right>=\left<u,v\right>_{H^{s}\left(\mathbb{R}^{N}\right)}
-\int_{\mathbb{R}^{N}}Q\left|u^{+}\right|^{q-1}vdx,\quad v\in H^{s}\left(\mathbb{R}^{N}\right).
$$
We claim that $E$ satisfies the Palais-Smale condition at level $c$ (in the sequel $\left(PS\right)_{c}$ condition) for any
\begin{align}\label{J}
c<c^{\ast}=\left(\frac{1}{2}-\frac{1}{q}\right)S_{s,q}^{\frac{q}{q-2}}\left(\mathbb{R}^{N}\right).
\end{align}
Indeed, for any sequence $\left\{u_{n}\right\}\subset H^{s}\left(\mathbb{R}^{N}\right)$ such that
\begin{align*}
\left\{
\begin{array}{lcl}
\displaystyle E\left(u_{n}\right)\rightarrow c; \\
\\
\displaystyle E'\left(u_{n}\right)\rightarrow 0,
\end{array}
\right.
\end{align*}
we have, by standard computations, that $u_{n}$ stays bounded in $H^{s}\left(\mathbb{R}^{N}\right)$. Passing if necessary to the subsequence, we assume
\begin{align*}
\left\{
\begin{array}{lcl}
\displaystyle u_{n}\rightharpoonup u,&&\text{in}~H^{s}\left(\mathbb{R}^{N}\right); \\
\\
\displaystyle u_{n}\rightarrow u,&&\text{in}~L_{loc}^{q}\left(\mathbb{R}^{N}\right); \\
\\
\displaystyle u_{n}\rightarrow u,&&\text{a.e.}~\text{on}~\mathbb{R}^{N}.
\end{array}
\right.
\end{align*}
Obviously $u$ satisfies 
\begin{align}\label{AB}
\left(-\Delta\right)^{s} u+u=Q\left|u^{+}\right|^{q-2}u^{+}.
\end{align}
We obtain $E\left(u\right)\geq0$. Let $v_{n}=u_{n}-u$. On the one hand, by Brezis-Lieb Lemma, there holds
\begin{align}\label{K}
E\left(v_{n}\right)=E\left(u_{n}\right)-E\left(u\right)+o\left(1\right)\leq c+o\left(1\right).
\end{align}
On the other hand, by $\left(Q_{3}\right)$ and
\begin{align*}
\left\{
\begin{array}{lcl}
\displaystyle \left<E\left(u'_{n}\right),u_{n}\right>\rightarrow0; \\
\\
\displaystyle \left<E'\left(u\right),u\right>=0,
\end{array}
\right.
\end{align*}
we get
\begin{align*}
\left\|v_{n}\right\|^{2}_{H^{s}\left(\mathbb{R}^{N}\right)}
\leq\left\|v_{n}\right\|^{q}_{L^{q}\left(\mathbb{R}^{N}\right)}+o\left(1\right).
\end{align*}
As $n\to\infty$, suppose that $\left\|v_{n}\right\|_{H^{s}\left(\mathbb{R}^{N}\right)}\neq o\left(1\right)$ in $H^{s}\left(\mathbb{R}^{N}\right)$, so that
$$
\left\|v_{n}\right\|_{H^{s}\left(\mathbb{R}^{N}\right)}\geq S_{s,q}^{\frac{q}{2\left(q-2\right)}}\left(\mathbb{R}^{N}\right)+o\left(1\right),
$$
and thus
\begin{align}\label{L}
E\left(v_{n}\right)
=\left(\frac{1}{2}-\frac{1}{q}\right)\left\|v_{n}\right\|^{2}_{H^{s}\left(\mathbb{R}^{N}\right)}+o\left(1\right)
\geq c^{\ast}+o\left(1\right).
\end{align}
The inequalities \eqref{K} and \eqref{L} contradicts \eqref{J}, therefore $v_{n}\rightarrow 0$ in $H^{s}\left(\mathbb{R}^{N}\right)$.

For any $c<c^{\ast}$, let $u_{s,q}$ be a positive minimizer of $S_{s,q}\left(\mathbb{R}^{N}\right)$. Since $E$ satisfies the mountain pass geometry, by the Ekeland variational principle there exists a sequence $\left\{u_{n}\right\}\subset H^{s}\left(\mathbb{R}^{N}\right)$ such that
\begin{align*}
\left\{
\begin{array}{lcl}
\displaystyle E\left(u_{n}\right)\rightarrow c_{0}; \\
\\
\displaystyle E'\left(u_{n}\right)\rightarrow 0,
\end{array}
\right.
\end{align*}
where
$$
c_{0}=\inf_{\gamma\in\Gamma}\max_{t\in\left[0,1\right]} E\left(\gamma\left(t\right)\right)
$$
and
$$
\Gamma:=\left\{\gamma\in\mathcal{C}\left(\left[0,1\right],H^{s}\left(\mathbb{R}^{N}\right)\right)
:\gamma\left(0\right)=0,\gamma\left(1\right)=e\right\}.
$$
Notice that by $\left(Q_{1}\right)$ and $\left(Q_{2}\right)$, we have
\begin{align*}
c\leq\max_{t\geq0}E\left(tu_{s,q}\right)<c^{\ast}.
\end{align*}
By the mountain pass theorem, there exists a nontrivial solution $u$ of \eqref{AB}. Let us multiply the equation \eqref{AB} by $u^{-}$ and integrate by parts, we have $u^{-}=0$, which means $u$ is a nonnegative solution of \eqref{A1}, and thus to be positive by the maximum principle for the fractional laplacian established by Cabr\'{e}-Sire \cite{Cabre-Sire}. By $E\left(u\right)<c^{\ast}$, we conclude that
$$
\left\|u\right\|^{2}_{H^{s}\left(\mathbb{R}^{N}\right)}<S_{s,q}^{\frac{q}{q-2}}\left(\mathbb{R}^{N}\right).
$$

Next we prove Theorem \ref{thm4}-$\left(2\right)$. Let us define the following energy functional $I\in \mathcal{C}^{1}\left(H^{s}\left(\mathbb{R}^{N}\right),\mathbb{R}\right)$ by
$$
I\left(u\right)=\frac{1}{2}\int_{\mathbb{R}^{N}}
\left|\left(-\Delta\right)^{\frac{s}{2}}u\right|^{2}+V\left|u\right|^{2}dx
$$
with
$$
\left<I'\left(u\right),v\right>
=\int_{\mathbb{R}^{N}}\left(-\Delta\right)^{\frac{s}{2}}u\left(-\Delta\right)^{\frac{s}{2}}v
+Vuvdx.
$$
We also define the unit sphere mainfold
$$
\mathcal{M}=\left\{u\in H^{s}\left(\mathbb{R}^{N}\right),\left\|u\right\|_{L^{q}\left(\mathbb{R}^{N}\right)}=1\right\}
$$
and
$$
I_{0}=\inf_{u\in\mathcal{M}}I\left(u\right).
$$
By $\left(V_{1}\right)$ and $\left(V_{2}\right)$, we have
\begin{align*}
I_{0}<\frac{1}{2}S_{s,q}\left(\mathbb{R}^{N}\right).
\end{align*}
Let $\left\{u_{n}\right\}\subset \mathcal{M}$ be a minimizing sequence of $I_{0}$, we obtain it is standard to prove that $\left\{u_{n}\right\}$ is bounded in $H^{s}\left(\mathbb{R}^{N}\right)$. Passing if necessary to the subsequence, we may assume
\begin{align*}
\left\{
\begin{array}{lcl}
\displaystyle u_{n}\rightharpoonup u,&&\text{in}~H^{s}\left(\mathbb{R}^{N}\right); \\
\\
\displaystyle u_{n}\rightarrow u,&&\text{in}~L_{loc}^{2}\left(\mathbb{R}^{N}\right); \\
\\
\displaystyle u_{n}\rightarrow u,&&\text{a.e.}~\text{on}~\mathbb{R}^{N}.
\end{array}
\right.
\end{align*}
Let $v_{n}=u_{n}-u$, by Brezis-Lieb Lemma, we get
\begin{align}\label{M}
1=\left\|u_{n}\right\|^{q}_{L^{q}\left(\mathbb{R}^{N}\right)}
=\left\|u\right\|^{q}_{L^{q}\left(\mathbb{R}^{N}\right)}+\left\|v_{n}\right\|^{q}_{L^{q}\left(\mathbb{R}^{N}\right)}+o\left(1\right),
\end{align}
moreover, we obtain
\begin{align}\label{N}
I\left(v_{n}\right)=I\left(u_{n}\right)-I\left(u\right)+o\left(1\right).
\end{align}
Notice that
\begin{align}\label{O}
I\left(u\right)\geq I_{0}\left\|u\right\|^{2}_{L^{q}\left(\mathbb{R}^{N}\right)}
\end{align}
and by $\left(V_{3}\right)$, we have
\begin{align}\label{P}
I\left(v_{n}\right)\geq\frac{1}{2}S_{s,q}\left(\mathbb{R}^{N}\right)\left\|v_{n}\right\|^{2}_{L^{q}\left(\mathbb{R}^{N}\right)}+o\left(1\right).
\end{align}
By joining \eqref{M}, \eqref{N}, \eqref{O} and \eqref{P} we conclude that
\begin{align*}
I_{0}=I\left(u_{n}\right)+o\left(1\right)\geq I_{0}\left\|u\right\|^{2}_{L^{q}\left(\mathbb{R}^{N}\right)}
+\frac{S_{s,q}}{2}\left(1-\left\|u\right\|^{q}_{L^{q}\left(\mathbb{R}^{N}\right)}\right)^{\frac{2}{q}}+o\left(1\right).
\end{align*}
This means that $\left\|u\right\|_{L^{q}\left(\mathbb{R}^{N}\right)}=1$ and then $u\in\mathcal{M}$. By the lower semicontinuity of $I$, we have that $u$ is a minimizer for $I$ in $\mathcal{M}$. Finally, assume $u$ is nonnegative, by the Lagrange multiplier rule, we have that
$$
u_{0}=(2I_{0})^{\frac{1}{q-2}}u
$$
is a positive solution of \eqref{A1} by the maximum principle again, which satisfies
$$
\left\|u_{0}\right\|_{L^{q}\left(\mathbb{R}^{N}\right)}<S_{s,q}^{\frac{1}{q-2}}\left(\mathbb{R}^{N}\right).
$$
\subsection{A nonlocal scalar field equation: proof of Theorem \ref{thm5}}
Define the functional $J\in \mathcal{C}^{1}\left(H^{\frac{1}{2}}\left(\mathbb{R}\right),\mathbb{R}\right)$ as follows
$$
J\left(u\right)=\frac{1}{2}\int_{\mathbb{R}}\left|\left(-\Delta\right)^{\frac{1}{4}}u\right|^{2}dx
$$
and the manifold
$$
\mathcal{N}=\left\{u\in H^{\frac{1}{2}}\left(\mathbb{R}\right)\backslash\left\{0\right\},
\frac{1}{2}\int_{\mathbb{R}}\left|u\right|^{2}dx
=\int_{\mathbb{R}}F\left(u\right)dx\right\},
$$
where $F\left(t\right)=\int^{t}_{0}f\left(s\right)ds$. Note that $\mathcal{N}$ is a $\mathcal{C}^{1}$ mainfold and nonempty by $\left(f_{1}\right)$. Let us denote
$$
J_{0}=\inf_{u\in \mathcal{N}}J\left(u\right)
$$
and
$$
U_{\frac{1}{2},q}=\frac{u_{\frac{1}{2},q}}{\left\|u_{\frac{1}{2},q}\right\|_{H^{\frac{1}{2}}\left(\mathbb{R}\right)}}\in H^{\frac{1}{2}}\left(\mathbb{R}\right),
$$
where $u_{\frac{1}{2},q}$ is a minimizer for $S_{\frac{1}{2},q}(\mathbb{R})$. By $\left(f_{3}\right)$, we have
\begin{align*}
\int_{\mathbb{R}}F\left(U_{\frac{1}{2},q}\right)dx\geq\frac{1}{2}>\frac{1}{2}\left\|U_{\frac{1}{2},q}\right\|^{2}_{L^{2}\left(\mathbb{R}\right)},
\end{align*}
and by $\left(f_{2}\right)$, for any sufficiently small $t_{1}>0$, we also have
\begin{align*}
\int_{\mathbb{R}}F\left(t_{1}U_{\frac{1}{2},q}\right)dx
<\frac{1}{2}\left\|t_{1}U_{\frac{1}{2},q}\right\|^{2}_{L^{2}\left(\mathbb{R}\right)}\ .
\end{align*}
Hence the continuous function
$$
h_{1}\left(t\right)=\frac{1}{2}\left\|tU_{\frac{1}{2},q}\right\|^{2}_{L^{2}
\left(\mathbb{R}\right)}-\int_{\mathbb{R}}F\left(tU_{\frac{1}{2},q}\right)dx
$$
satisfies
\begin{align*}
\left\{
\begin{array}{lcl}
\displaystyle h_{1}\left(0\right)=0; \\
\\
\displaystyle h_{1}\left(t_{1}\right)>0;\\
\\
\displaystyle h_{1}\left(1\right)<0,
\end{array}
\right.
\end{align*}
which means that there exists $\overline{t_{1}}\in\left(0,1\right)$ such that $h_{1}\left(\overline{t_{1}}\right)=0$ and then $\overline{t_{1}}U_{\frac{1}{2},q}\in \mathcal{N}$, therefore we have
\begin{align*}
J_{0}\leq J\left(\overline{t_{1}}U_{\frac{1}{2},q}\right)<\frac{1}{2}.
\end{align*}
Let $\left\{u_{n}\right\}\subset \mathcal{N}$ be a minimizing sequence for $J_{0}$ and consider the radially symmetric sequence
$$
U_{n}=u_{n}^{\#}\left(\lambda \left\|u_{n}^{\#}\right\|^{2}_{L^{2}\left(\mathbb{R}\right)}x\right),
$$
where
$$
\lambda\geq\frac{1}{1-2J_{0}},
$$
to obtain $\left\{U_{n}\right\}\subset \mathcal{N}$ and
\begin{align*}
\lim_{n\rightarrow+\infty}J\left(U_{n}\right)\leq J_{0}\ .
\end{align*}
Therefore $\left\{U_{n}\right\}$ is also a minimizing sequence for $J_{0}$. Passing if necessary to a subsequence, we have
\begin{align*}
\sup_{n}\left\|U_{n}\right\|^{2}_{H^{\frac{1}{2}}\left(\mathbb{R}\right)}
\leq\sup_{n}\left(2J_{0}+\frac{1}{\lambda}+o\left(1\right)\right)\leq1
\end{align*}
and we may assume \begin{align*}
\left\{
\begin{array}{lcl}
\displaystyle U_{n}\rightharpoonup U,&&\text{in}~H_{rad}^{\frac{1}{2}}\left(\mathbb{R}\right); \\
\\
\displaystyle U_{n}\rightarrow U,&&\text{a.e.}~\text{on}~\mathbb{R}.
\end{array}
\right.
\end{align*}
Let
\begin{align*}
Q\left(t\right)=e^{\pi t^{2}}-1,
\end{align*}
notice that by $\left(f_{1}\right)$, $\left(f_{2}\right)$ and L'H\^{o}pital's rule, we have
\begin{align}\label{R}
\lim_{t\rightarrow0}\frac{F\left(t\right)}{Q\left(t\right)}=\lim_{\left|t\right|\rightarrow+\infty}\frac{F\left(t\right)}{Q\left(t\right)}=0,
\end{align}
and by Proposition \ref{Prop2}, we know
\begin{align}\label{S}
\sup_{n}\int_{\mathbb{R}}\left|Q\left(U_{n}\right)\right|dx<+\infty,
\end{align}
therefore \eqref{R} and \eqref{S} together with
\begin{align*}
F\left(U_{n}\right)\rightarrow F\left(U\right)\quad \text{a.e.}~\text{on}~\mathbb{R}
\end{align*}
satisfy the conditions of the compactness lemma of Strauss established by Berestycki-Lions \cite{Berestycki-Lions}, which gives
$$
\int_{\mathbb{R}}F\left(U_{n}\right)dx\rightarrow \int_{\mathbb{R}}F\left(U\right)dx,
$$
and in turn we have
\begin{align*}
\frac{1}{2}\left\|U\right\|_{L^{2}\left(\mathbb{R}\right)}
\leq\frac{1}{2}\liminf_{n\rightarrow+\infty}\left\|U_{n}\right\|_{L^{2}\left(\mathbb{R}\right)}
=\int_{\mathbb{R}}F\left(U\right)dx.
\end{align*}
If the strict inequality holds
$$
\frac{1}{2}\left\|U\right\|_{L^{2}\left(\mathbb{R}\right)}<\int_{\mathbb{R}}F\left(U\right)dx,
$$
by $\left(f_{1}\right)$ for sufficiently small $t_{2}>0$, we get
\begin{align*}
\int_{\mathbb{R}}F\left(t_{2}U\right)dx<\frac{1}{2}\left\|t_{2}U\right\|^{2}_{L^{2}\left(\mathbb{R}\right)},
\end{align*}
so that the continuous function
$$
h_{2}\left(t\right)=\frac{1}{2}\left\|tU\right\|^{2}_{L^{2}\left(\mathbb{R}\right)}-\int_{\mathbb{R}}F\left(tU\right)dx
$$
satisfies
\begin{align*}
\left\{
\begin{array}{lcl}
\displaystyle h_{2}\left(0\right)=0; \\
\\
\displaystyle h_{2}\left(t_{2}\right)>0;\\
\\
\displaystyle h_{2}\left(1\right)<0,
\end{array}
\right.
\end{align*}
hence there exists $\overline{t_{2}}\in\left(0,1\right)$ such that $h_{2}\left(\overline{t_{2}}\right)=0$ and then $\overline{t_{2}}U\in \mathcal{N}$, which leads to
\begin{align*}
J_{0}&\leq J\left(t_{2}U\right)<\liminf_{n\rightarrow+\infty}J\left(U_{n}\right)=J_{0},
\end{align*}
and then a contradiction. Therefore we deduce that
$$
\frac{1}{2}\left\|U\right\|_{L^{2}\left(\mathbb{R}\right)}=\int_{\mathbb{R}}F\left(U\right)dx,
$$
that is $U\in \mathcal{N}$. By lower semicontinuity of $J$, $U$ is a minimizer for $J$ in $\mathcal{N}$. Following the standard preceedure of Berestycki-Lions \cite{Berestycki-Lions}, we also have that $U$ is a positive ground state solution of \eqref{A2}.

\subsection{A nonlocal Schr\"{o}dinger system: proof of Theorem \ref{thm7}}
Let us recall from Lu-Peng \cite{Lu-Peng} the following system
\begin{align}\label{B3}
\left\{
\begin{array}{lcl}
\displaystyle \left(-\Delta\right)^{s} u+u=f\left(u\right)+\lambda v,\\
\\
\displaystyle \left(-\Delta\right)^{s} v+v=g\left(v\right)+\lambda u,\\
\\
\displaystyle u,v\in H^{s}\left(\mathbb{R}^{N}\right),
\end{array}
\right.
\end{align}
where $0<s<1$, $f,g\in\mathcal{C}^{1}\left(\mathbb{R}\right)$ and $0<\lambda<1$, for which they proved the following
\begin{Prop}\label{Prop3}
Let $f,g$ satisfy
\begin{itemize}
\item $\displaystyle \lim_{t\rightarrow0^{+}}\frac{f\left(t\right)}{t}=\lim_{t\rightarrow0^{+}}\frac{g\left(t\right)}{t}=0$;
\item there exist $p,q\in\left(2,2_{s}^{\ast}\right)$ such that
$$
\lim_{\left|t\right|\rightarrow+\infty}\frac{f\left(t\right)}{\left|t\right|^{p-1}}
=\lim_{\left|t\right|\rightarrow+\infty}\frac{g\left(t\right)}{\left|t\right|^{q-1}}
=0;
$$
\item there exist $\zeta_{1},\zeta_{2}>0$ such that
\begin{align*}
\left\{
\begin{array}{lcl}
\displaystyle \int^{\zeta_{1}}_{0}f\left(t\right)dt>\frac{\zeta_{1}^{2}}{2},\\
\\
\displaystyle \int^{\zeta_{2}}_{0}g\left(t\right)dt>\frac{\zeta_{2}^{2}}{2}.
\end{array}
\right.
\end{align*}
\end{itemize}
For any $\lambda\in\left(0,1\right)$, the system \eqref{B3} has a ground state solution.
\end{Prop}
Let
\begin{align*}
\left\{
\begin{array}{lcl}
\displaystyle f\left(t\right)=\left|t\right|^{p-2}t;\\
\\
\displaystyle g\left(t\right)=\left|t\right|^{q-2}t,
\end{array}
\right.
\end{align*}
then the system \eqref{A4} satisfies Proposition \ref{Prop3} and thus we have Theorem \ref{thm7}-$\left(1\right)$. Next we prove Theorem \ref{thm7}-$\left(2\right)$. The authors in \cite{Zhen-He-Xu-Yang} studied the following Schr\"{o}dinger system
\begin{align}\label{B4}
\left\{
\begin{array}{lcl}
\displaystyle \left(-\Delta\right)^{s} u+\alpha u=\left|u\right|^{2^{\ast}_{s}-2}u+\lambda v,\\
\\
\displaystyle \left(-\Delta\right)^{s} v+\beta v=\left|v\right|^{q-2}v+\lambda u,\\
\\
\displaystyle u,v\in H^{s}\left(\mathbb{R}^{N}\right),
\end{array}
\right.
\end{align}
and obtaining teh following result in the spirit of Theorem \ref{thm6}.
\begin{Prop}\label{Prop4}
Let $0<s<1$, $N>2s$, $\alpha,\beta>0$ and $0<\lambda<\sqrt{\alpha\beta}$. Let
$$
\alpha_{s}=\left[\frac{\mathcal{S}_{s}^{\frac{N}{2s}}}{\frac{N}{s}\left(\frac{1}{2}-\frac{1}{q}\right)
S_{s,q}^{\frac{q}{q-2}}\left(\mathbb{R}^{N}\right)}\right]^{\frac{1}{\frac{q}{q-2}-\frac{N}{2s}}},
$$
\item [$\left(1\right)$] if $0<\alpha\leq\alpha_{s}$, then \eqref{B4} has a positive radial decreasing ground state solution.
\item [$\left(2\right)$] if $\alpha>\alpha_{s}$, then there exists $\lambda_{s}\in\left[\sqrt{\left(\alpha-\alpha_{0}\right)\beta},\sqrt{\alpha\beta}\right)$ such that
\begin{itemize}
\item [$\left(I\right)$] if $\lambda<\lambda_{s}$, then \eqref{B4} has no ground state solution;
\item [$\left(II\right)$] if $\lambda>\lambda_{s}$, then \eqref{B4} has a positive radial decreasing ground state solution.
\end{itemize}
\end{Prop}

Let
\begin{align*}
\left\{
\begin{array}{lcl}
\displaystyle \tau_{1}=\frac{N}{s}\left(\frac{1}{2}-\frac{1}{q}\right),\\
\\
\displaystyle \tau_{2}=\frac{N}{s}\left(\frac{1}{q}-\frac{1}{2_{s}^{\ast}}\right),
\end{array}
\right.
\end{align*}
for which the following holds
$$
\tau_{1}+\tau_{2}=1.
$$
By Theorem \ref{thm2}-$\left(3\right)$, we have
\begin{align*}
S_{s,q}\left(\mathbb{R}^{N}\right)\geq\tau_{2}^{-\tau_{2}}\tau_{1}^{-\tau_{1}}\mathcal{S}_{s}^{\tau_{1}},
\end{align*}
hence
\begin{align*}
\alpha_{s}\leq\tau_{2}\tau_{1}^{\frac{\frac{N}{2s}-1}{\frac{q}{q-2}-\frac{N}{2s}}}<1\ .
\end{align*}
As a consequence, \eqref{A4} satisfies Proposition \ref{Prop4}-$\left(2\right)$ and this proves Theorem \ref{thm7}-$\left(2\right)$. Finally, we prove Theorem \ref{thm7}-$\left(3\right)$. Let us use as test functions in  \eqref{A4}, $u$ and $v$ to get
\begin{align}\label{T}
\left\{
\begin{array}{lcl}
\displaystyle \int_{\mathbb{R}^{N}}\left|\left(-\Delta\right)^{\frac{s}{2}}u\right|^{2}dx
+\int_{\mathbb{R}^{N}}\left|u\right|^{2}dx
=\int_{\mathbb{R}^{N}}\left|u\right|^{2^{\ast}_{s}}dx
+\lambda\int_{\mathbb{R}^{N}}uvdx,\\
\\
\displaystyle \int_{\mathbb{R}^{N}}\left|\left(-\Delta\right)^{\frac{s}{2}}v\right|^{2}dx
+\int_{\mathbb{R}^{N}}\left|v\right|^{2}dx
=\int_{\mathbb{R}^{N}}\left|v\right|^{2^{\ast}_{s}}dx
+\lambda\int_{\mathbb{R}^{N}}uvdx.
\end{array}
\right.
\end{align}
Moreover, by the Poho\v{z}aev identity established by Chang-Wang \cite{Chang-Wang}, we also have
\begin{align}\label{U}
\int_{\mathbb{R}^{N}}\left|\left(-\Delta\right)^{\frac{s}{2}}u\right|^{2}dx
&+\int_{\mathbb{R}^{N}}\left|\left(-\Delta\right)^{\frac{s}{2}}v\right|^{2}dx
+\frac{2^{\ast}_{s}}{2}\left(\int_{\mathbb{R}^{N}}\left|u\right|^{2}dx
+\int_{\mathbb{R}^{N}}\left|v\right|^{2}dx\right)\nonumber\\
&=\int_{\mathbb{R}^{N}}\left|u\right|^{2^{\ast}_{s}}dx
+\int_{\mathbb{R}^{N}}\left|v\right|^{2^{\ast}_{s}}dx
+2^{\ast}_{s}\lambda\int_{\mathbb{R}^{N}}uvdx.
\end{align}
By combining \eqref{T} and \eqref{U} we obtain
$$
\int_{\mathbb{R}^{N}}\left|u\right|^{2}dx
+\int_{\mathbb{R}^{N}}\left|v\right|^{2}dx=
2\lambda\int_{\mathbb{R}^{N}}uvdx,
$$
which yields $u=v=0$ by Cauchy-Schwarz inequality and the assumption $0<\lambda<1$.

\end{document}